\documentclass[review]{elsarticle}

\usepackage{lineno,hyperref}
 \usepackage{amsmath, amsthm, amsfonts}
 \usepackage{mathrsfs, amssymb}
 \usepackage{graphicx, color, bm, cases}
\modulolinenumbers[5]

\newtheorem{thm}{Theorem}[section]

\newtheorem{lemma}[thm]{Lemma}

\newtheorem{prop}[thm]{Proposition}
\newtheorem{rem}[thm]{Remark}

\newtheorem{con}[thm]{Condition}
\renewcommand
{\theequation}{\thesection.\arabic{equation}}

\def\ga{{\gamma}}

\def\<{\left<}\def\>{\right>}\def\({\left(}\def\){\right)}

\numberwithin{equation}{section}
\newfam\msbmfam\font\tenmsbm=msbm10\textfont
\msbmfam=\tenmsbm\font\sevenmsbm=msbm7

\scriptfont\msbmfam=\sevenmsbm

 \def\beqlb{\begin{eqnarray}}\def\eeqlb{\end{eqnarray}}
 \def\beqnn{\begin{eqnarray*}}\def\eeqnn{\end{eqnarray*}}
 \def\d{{\mbox{\rm d}}}\def\e{{\mbox{\rm e}}}

 \def\mcr{\mathscr}\def\mbb{\mathbb}
\def\goto{{\rightarrow}}
 \def\ar{\!\!&}\def\nnm{\nonumber}\def\ccr{\nnm\\}

\journal{Stochastic Processes and their Applications}

\begin{document}

\begin{frontmatter}

\title{{  Well-posedness of the martingale problem for super-Brownian motion with interactive branching}}


\author[mymainaddress]{Lina Ji}
\ead{jiln@smbu.edu.cn}

\author[mysecondaryaddress]{Jie Xiong}
\ead{xiongj@sustech.edu.cn}

\author[mythirdaryaddress]{Xu Yang\corref{mycorrespondingauthor}}
\cortext[mycorrespondingauthor]{Corresponding author}
\ead{xuyang@mail.bnu.edu.cn}

\address[mymainaddress]{{Faculty of Computational Mathematics and Cybernetics, Shenzhen MSU-BIT University, Shenzhen, China.}}
\address[mysecondaryaddress]{Department of Mathematics and National Center for Applied Mathematics (Shenzhen), Southern University of Science and Technology, Shenzhen, China.}
\address[mythirdaryaddress]{School of Mathematics and Information Science, North Minzu University, Yinchuan, China.}

\begin{abstract}
In this paper a martingale problem for super-Brownian motion with interactive branching is derived. The uniqueness of the solution to the martingale problem is obtained by using the pathwise uniqueness of the solution to a corresponding system of SPDEs with proper boundary conditions. The existence of the solution to the martingale problem and the local H\"{o}lder continuity of the density process are also studied.
\end{abstract}


\end{frontmatter}

\linenumbers

\section{Introduction and main results}
\setcounter{equation}{0}
\renewcommand{\theequation}{\thesection.\arabic{equation}}

Let $M(\mbb{R})$ be the collection of all finite Borel measures on $\mbb{R}$ endowed with the weak convergence topology. Let $C_b^n(\mbb{R})$ be the set of all bounded continuous functions on $\mbb{R}$ with continuous bounded derivatives up to $n$th order.
We consider a continuous $M(\mbb{R})$-valued process $(X_t)_{t \ge 0}$ satisfying the martingale problem (MP): for $\phi \in C_b^2(\mbb{R}),$ the process
\beqlb\label{0.1}
M_t(\phi)\equiv\<X_t,\phi\> - \<X_0, \phi\>-\int^t_0\<X_s,\frac12\phi''\>ds
\eeqlb
is a continuous martingale with quadratic variation process $(\<M(\phi)\>_t)_{t \ge 0}$ given by
\beqlb\label{0.2}
\<M(\phi)\>_t=\int_0^t\<X_s,\ga(\mu_s,\cdot)\phi^2\>ds,
\eeqlb
where $(\mu_t)_{t \ge 0}$ is the density process of $(X_t)_{t\ge 0},$ and $\gamma$ is the interacting branching rate depending on $(\mu_t)_{t \ge 0}$. The notation $\< \nu, \phi\>$ denotes the integral of the function $\phi$ with respect to the measure $\nu.$ In this paper we always assume that $X_0(d x) = \mu_0(x)d x$ with $\mu_0 \in C_c(\mbb{R})^+,$ where $C_c(\mbb{R})^+$ is the collection of nonnegative continuous functions on $\mbb{R}$ with compact support. Obviously, $\mu_0 \in L^1(\mbb{R})^+$ with $L^1(\mbb{R})^+$ being the set of nonnegative functions $f$ on $\mbb{R}$ with $\int_{\mbb{R}}f(x)d x < \infty.$ 
When $\gamma$ is a constant,  the process $(X_t)_{t \ge 0}$ is a classical super-Brownian motion. In this case the well-posedness of the MP (\ref{0.1}, \ref{0.2}) was established using the nonlinear partial differential equation satisfied by its log-Laplace transform, which was obtained by Watanabe \cite{W68}. Moreover, a new approach for the well-posedness of the MP was suggested by Xiong~\cite{X13b},  in which a relationship between the classical super-Brownian motion and a  stochastic partial differential equation (SPDE) satisfied by its corresponding distribution-function-valued process was established. The weak uniqueness of the solution to the MP was also obtained by the strong uniqueness of the solution to the corresponding SPDE in  \cite{X13b}. See He et al.~\cite{HLY14} for the case of super-L\'{e}vy process.


 Superprocesses with interaction are natural to consider, since for many species branching and spatial motion depend on the population density. When the spatial motion is interactive,  the well-posedness of the martingale problem was studied by Donnelley and Kurtz~\cite{DK99}, see also 
 Perkins~\cite[Theorem~V.5.1]{P95} and Li et~al. \cite{LWX05}. Uniqueness for the historical superprocesses with certain interaction was investigated by Perkins~\cite{P95}. Furthermore, the superprocesses with interactive immigration was studied in \cite{ DL03, FL04,S90}, see also Li~\cite[Section~10]{L11}. The well-posedness of the martingale problem for the interactive immigration process was solved by Mytnik and Xiong \cite{MX15}. See also \cite{XY16} for the well-posedness of the martingale problem for a superprocess with location-dependent branching, interactive immigration mechanism and spatial motion.

However, the hard case for the superprocess with interactive branching was rarely investigated. We are interested in the case of $\gamma(\mu_s, x) = \gamma(\mu_s(x)),$ i.e., the interactive branching mechanism depending on the density process $(\mu_t)_{t \ge 0},$ where the well-posedness of the MP (\ref{0.1}, \ref{0.2}) is still an {\it open problem}. The weak uniqueness of the solution to the MP (\ref{0.1}, \ref{0.2}) is very difficult to prove. As a first step, throughout this paper we assume that $\gamma$ satisfies the following condition:

 \begin{con}\label{0.4b}
Let integer $n \ge 0$ and $-\infty=a_0<\cdots<a_{n+1}=\infty$ be fixed. For any $f \in L^1(\mbb{R})^+$, let
 \beqnn
\ga(f,x)=
\begin{cases}
g_i^2(f(a_{i + 1})),& a_i\le x < a_{i+1},\ i = 0, \cdots, n - 1,\\
 g_n^2\(\int_{a_n}^\infty f(y)d y\), & a_n \le x < \infty,
\end{cases}
\eeqnn
where $g_i,\ i = 0, \cdots, n$ are continuous bounded functions from $\mbb{R}_+$ to $\mbb{R}_+.$
\end{con}
The existence and local H\"{o}lder continuity of the density process $(\mu_t)_{t \ge 0}$ are investigated in this paper, where $(\mu_t)_{t \ge 0}$ is a weak solution to the following SPDE:
\beqlb\label{1.1}
\frac{\partial}{\partial t} \mu_t(x) = \frac{1}{2}\Delta \mu_t(x) + \sqrt{\mu_t(x)\gamma(\mu_t, x)}\dot{W}_{tx},
\eeqlb
$\Delta$ denotes the one-dimensional Laplacian operator and $\dot{W}_{tx}$ is the formal derivative of a white noise random measure. Let $C(\mbb{R})^+$ be the collection of nonnegative continuous functions on $\mbb{R}$. The notation $\langle f, g\rangle$ denotes $\int_{\mbb{R}} f(x)g(x) d x$  if it exists. The $C(\mbb{R})^+$-valued process $(\mu_t)_{t \ge 0}$ is a weak solution of the SPDE \eqref{1.1}  in the following sense: for any $\phi \in C_b^2(\mbb{R}),$ we have
\beqlb\label{1.1+}
\langle \mu_t, \phi\rangle \ar=\ar \langle \mu_0, \phi\rangle + \frac{1}{2}\int_0^t \langle \mu_s, \phi''\rangle d s\ccr
\ar\ar + \int_0^t\int_{\mbb{R}} \sqrt{\mu_s(x)\gamma(\mu_s, x)}\phi(x)W(d s, d x), \quad t \ge 0
\eeqlb
almost surely, where $W(d s, d x)$ is a time-space Gaussian white noise random measure with the Lebesgue measure $d s d x$ as its intensity. The existence of the solution to the MP (\ref{0.1}, \ref{0.2}) is given by showing the existence of the solution to \eqref{1.1}. Moreover, we show the weak uniqueness of the solution to the MP~(\ref{0.1}, \ref{0.2}) in Theorem~\ref{main} below. The main idea is to relate the MP to a system of SPDEs, which is satisfied by a sequence of corresponding distribution-function-valued processes on intervals. The weak uniqueness of solution to the MP follows from the pathwise uniqueness of the solution to the system of SPDEs, see Section 3. Throughout this paper we always assume that all random variables defined on the same filtered probability space $(\Omega, \mcr{F}, \mcr{F}_t, \mbb{P})$ without special explanation. Let $\mbb{E}$ be the corresponding expectation.


\begin{thm}\label{t1}
There exists a continuous $M(\mbb{R})$-valued process $(X_t)_{t \ge 0}$ satisfying the MP (\ref{0.1}, \ref{0.2}), where $X_t(d x)$ is absolutely continuous respect to $d x$ with density $\mu_t(x)$ satisfying \eqref{1.1} for all $t \ge 0$ almost surely, and
\beqlb\label{moment}
 \sup_{0 \le t \le T}\mbb{E}\left[\int_{\mbb{R}}\mu_t(x)^{2p}e^{-|x|}d x\right] < \infty
 \eeqlb
 for every $T > 0$ and $p \ge 1$.
\end{thm}

\begin{thm}\label{t2}
Suppose that $(\mu_t)_{t \ge 0}$ satisfies \eqref{1.1} and \eqref{moment}. Then $(\mu_t)_{t > 0}$ is locally H{\"o}lder continuous with exponent $\lambda_1/2$ in time variable and with exponent $\lambda_2$ in space variable for all $\lambda_1, \lambda_2 \in (0, 1/2).$ Namely,  for any fixed $0< r_0 < T$ and $L > 0,$  there exists a random variable $\xi \ge 0$ depending on $\lambda_1, \lambda_2$ such that with probability one,
\beqnn
|\mu_t(x) - \mu_r(y)| \le \xi (|t - r|^{\lambda_1/2} + |x - y|^{\lambda_2}), \quad t, r \in [r_0, T],\ x, y \in [-L, L].
\eeqnn
\end{thm}

\begin{thm}\label{main}
 Assume that $g_n$ in Condition \ref{0.4b} is $\beta$-H\"{o}lder continuous with $ \frac{1}{2}\le \beta \le 1$, i.e.,
\beqlb\label{beta}
|g_n(x) - g_n(y)| \le K |x - y|^\beta, \qquad  x, y \ge 0
\eeqlb
for some constant $K > 0.$ Then the weak uniqueness of the solution to the MP (\ref{0.1}, \ref{0.2}) holds.
\end{thm}

\begin{rem}
Taking $n = 0$ in Condition~\ref{0.4b}, the branching rate depends on the total mass process, i.e., the quadratic variation process $(\langle M(\phi)\rangle_t)_{t \ge 0}$ of the martingale defined by \eqref{0.1} is
\beqlb\label{rem}
\langle M(\phi)\rangle_t = \int_0^t \gamma(\langle X_s, 1\rangle )\langle X_s, \phi^2\rangle d s.
\eeqlb
The well-posedness of the MP (\ref{0.1}, \ref{rem}) is a corollary of Theorems~\ref{t1} and \ref{main}.
\end{rem}

We now introduce some notation. 
 Let $B(\mbb{R})$ be the collection of all bounded functions on $\mbb{R}.$ We use $B(\mbb{R})^+$ to denote the subset of $B(\mbb{R})$ of nonnegative functions.
Let $B[a_1, a_2]$ be the Banach space of bounded measurable functions on $[a_1, a_2]$ furnished with the supremum norm. For $f, g \in B [a_1, a_2]$ let $\< f, g\> = \int_{a_1}^{a_2} f(x)g(x) d x$. Define $C_b[a_1, a_2]$ to be the set of bounded continuous functions on $[a_1, a_2]$. For any integer $n \ge 0,$ let $C_b^n[a_1, a_2]$ be the subset of $C_b[a_1, a_2]$ of functions with bounded continuous derivatives up to the $n$th order. Let $C_c^n(a_1, a_2)$ denote the subset of $C_b^n[a_1, a_2]$ of functions with compact support in $(a_1, a_2).$ Define $C_c^\infty(\mbb{R})$ to be the set of infinitely differentiable functions on $\mbb{R}$ with compact support.  For $f, g \in C(\mbb{R})^+$ and $\lambda > 0,$ let $|f - g|_{(-\lambda)} = \sup_{x \in \mbb{R}}|e^{-\lambda |x|}(f(x) - g(x))|.$ Let $C_{\rm tem}^+(\mbb{R})$ be the subspace of $C(\mbb{R})^+$ of functions $f$ with $|f|_{(-\lambda)}< \infty$ for every $\lambda > 0$, whose topology is induced by norm $\{|f|_{(-\lambda)}: \lambda > 0\}.$ Let  $C([0, T]\times [-L, L], \mbb{R}_+)$ be the space of continuous functions from $[0, T]\times [-L, L]$ to $\mbb{R}_+$ furnished with the supremum norm $\|\cdot\|_{[0,T]\times [-L, L]}$, and $C([0, T]\times \mbb{R}, \mbb{R}_+)$ be the space of continuous functions from $[0, T]\times \mbb{R}$ to $\mbb{R}_+$ furnished with metric
\beqlb\label{dfg}
d(f,g) = \int_0^\infty e^{-L} \left(\|f-g\|_{[0,T]\times [-L, L]} \wedge 1\right) d L.
\eeqlb
For a Banach space $\mathcal{Y}$ with suitable norm $\|\cdot\|_\mathcal{Y}$, let $C([0, T],  \mathcal{Y})$ 
be the set of all continuous maps from $[0, T]$ to $\mathcal{Y}$
with the topology induced by 
\beqlb\label{d1fg}
d_1(f,g) =  \sup_{0 \le t \le T}\|f(t)-g(t)\|_\mathcal{Y}.
\eeqlb

The rest of the paper is organized as follows. In Section~2, we present the proofs of Theorems~\ref{t1} and \ref{t2}. The weak uniqueness of the solution to the MP~(\ref{0.1}, \ref{0.2}), i.e., the proof of Theorem~\ref{main}, is given in Section~3. Moreover, several useful lemmas and proofs are given in Section 4. Throughout the paper we use $\nabla$ to denote the first order spatial differential operator. We use $K$ to denote a non-negative constant whose value may change from line to line. In the integrals, we make convention that, for $a \le b \in \mbb{R},$
\beqnn
\int_a^b = \int_{(a, b]}\quad \text{and}\quad \int_a^\infty = \int_{(a, \infty)}.
\eeqnn

\section{Proofs of Theorems~\ref{t1} and \ref{t2}}

In this section, we present the proofs of Theorems~\ref{t1} and \ref{t2}. The existence of the solution to the MP (\ref{0.1}, \ref{0.2}) is obtained by the existence of the corresponding density process. Moreover, the local H{\"o}lder continuity of the density process $(\mu_t)_{t \ge 0}$ is given. 
\begin{lemma}\label{t2.3}
The martingale defined in \eqref{0.1} induces an $(\mcr{F}_t)$-martingale measure  satisfying
\beqnn
M_t(\phi) = \int_0^t \int_{\mbb{R}} \phi(x) M(d s, d x), \qquad t \ge 0,\ \phi \in C_b^2(\mbb{R}),
\eeqnn
where $M(d s, d x)$ is an orthogonal martingale measure on $\mbb{R}_+\times\mbb{R}$ with covariance measure $d s\int_{\mbb{R}}\big[\gamma(\mu_s, z) \delta_z(d x)\delta_z(d y)\big]X_s(d z)$.
\end{lemma}

 The definition of orthogonal martingale measure mentioned above is given in, e.g.,  Walsh \cite[p. 288]{W86}.

\noindent{\it Proof of Lemma~\ref{t2.3}.}  By the MP~(\ref{0.1}, \ref{0.2}), one can see that
\beqlb\label{Gamma1}
\mbb{E}[\<X_t, 1\>] = \<X_0, 1\> < \infty.
\eeqlb
For each $n \ge 1$ we define the measure $\Gamma_n \in M(\mbb{R})$ by
\beqlb\label{Gamma2}
\Gamma_n(\phi) := \mbb{E}\Big[\int_0^n d s \int_{\mbb{R}}\gamma(\mu_s, x) \phi(x)X_s(d x)\Big]
\eeqlb
with $\phi \in B(\mbb{R})^+.$ Note that $\gamma$ is bounded by Condition \ref{0.4b}. It then follows from \eqref{Gamma1} and \eqref{Gamma2} that for each $n \ge 1,$ $\Gamma_n(\phi) \le K \int_0^n \mbb{E} [\<X_s, 1\>] d s \le Kn\<X_0, 1\>$ is bounded. The rest proof follows by changing $c(z) = \gamma(\mu_s, z)/2$ and $H(z, d \nu) = 0$ in the proof of Li~\cite[Theorem~7.25]{L11}. We omit it here.
\qed

\begin{lemma}\label{<mu,1>}
 Suppose that $(\mu_t)_{t \ge 0}$ is a solution to \eqref{1.1}. Then for every $t \ge 0,$ we have $\mbb{E}[\langle\mu_t, 1\rangle] = \langle \mu_0, 1\rangle.$
 \end{lemma}

 \proof
 By \eqref{1.1} we have
\beqnn
\langle \mu_t, 1\rangle = \langle \mu_0, 1\rangle + \int_0^t\int_{\mbb{R}} \sqrt{\mu_s(x)\gamma(\mu_s, x)} W(d s, d x).
\eeqnn
For each $n \ge 1$ we define stopping time $\tau_n$ by $\tau_n = \inf\{t \ge 0: \langle \mu_t, 1\rangle \ge n\}.$ For $T > 0,$ it follows from the continuity of $t \mapsto \langle \mu_t, 1\rangle$ that $\sup_{0 \le t \le T} \langle \mu_t, 1\rangle < \infty$ almost surely, which implies $\tau_n \rightarrow \infty$ almost surely as $n \rightarrow \infty.$ Further, $t \mapsto \int_0^{t\wedge \tau_n}\int_{\mbb{R}} \sqrt{\mu_s(x)\gamma(\mu_s, x)} W(d s, d x)$
is a martingale (see, e.g., \cite[p. 55]{IW89}), since
\beqnn
\mbb{E}\left[\int_0^{t\wedge \tau_n}\int_{\mbb{R}} \mu_s(x)\gamma(\mu_s, x) d s d x\right] \le K\mbb{E}\left[\int_0^{t\wedge \tau_n} \langle \mu_s, 1\rangle d s\right] \le Ktn.
\eeqnn
It then implies that $\mbb{E}\left[\langle \mu_{t\wedge \tau_n}, 1\rangle\right] = \langle \mu_0, 1\rangle.$ Further,
\beqnn
\mbb{E}[\langle \mu_{t\wedge \tau_n}, 1\rangle^2] \ar\le\ar 2\langle \mu_0, 1\rangle^2 + 2\mbb{E}\left[\left(\int_0^{t\wedge \tau_n}\int_{\mbb{R}} \sqrt{\mu_s(x)\gamma(\mu_s, x)} W(d s, d x)\right)^2\right]\cr
\ar\le\ar 2\langle \mu_0, 1\rangle^2 + K \int_0^t \mbb{E}[\langle \mu_{s\wedge \tau_n}, 1\rangle] d s\cr
\ar\le\ar K [\langle \mu_0, 1\rangle^2 + t\langle \mu_0, 1\rangle],
\eeqnn
which implies that $\sup_n \mbb{E}[\langle \mu_{t\wedge \tau_n}, 1\rangle^2] < \infty.$ By \cite[p. 67]{K02}, we have $\mbb{E}[\langle\mu_t, 1\rangle] = \lim_{n \rightarrow \infty}\mbb{E}\left[\langle \mu_{t\wedge \tau_n}, 1\rangle\right] = \langle \mu_0, 1\rangle.$ The proof ends here.
 \qed

 Let $P_t(x, d y)$ be the semigroup generated by $\Delta/2,$ which is absolutely continuous with respect to Lebesgue measure $d y$ with density $p_t(x, y)$ given by 
\[p_t(x, y) := p_t(x - y) :=\frac{1}{\sqrt{2\pi t}}\e^{-|x - y|^2/(2t)}, \qquad t > 0,\  x,y \in \mbb{R}.\]
Observe that
\beqlb\label{p_t}
p_t(x - y) \le \frac{1}{\sqrt{t}}, \qquad  t > 0,\ x,\ y \in \mbb{R}.
\eeqlb

 \begin{lemma}\label{0405}
Let $(X_t)_{t \ge 0}$ satisfy the MP (\ref{0.1}, \ref{0.2}) and $T > 0$ be fixed. Then for any $p \ge 1$ we have $
\mbb{E}\left[\sup_{0 \le t \le T}\<X_t, 1\>^{2p}\right] < \infty.$
\end{lemma}
\proof
For each $n \ge 1,$ we define stopping time $\kappa_n = \inf\{t \ge 0: \<X_t, 1\>^{2p} \ge n\}.$
Note that $\sup_{0 \le t \le T}\<X_t, 1\>^{2p} < \infty$ almost surely because of the continuity of $(\<X_t, 1\>)_{t \ge 0}$. We then have $\lim_{n \rightarrow \infty}(T \wedge \kappa_n) = T$ almost surely. By the MP (\ref{0.1}, \ref{0.2}) we have
\beqlb\label{0413a}
\<X_t, 1\> = \<X_0, 1\> + M_t(1).
\eeqlb
By Burkholder-Davis-Gundy's inequality, H{\"o}lder's inequality and $a \le a^2 + 1$ for $a \in \mbb{R},$ one can see that
\beqnn
\mbb{E}\left[\sup_{0 \le t \le T\wedge \kappa_n}|M_t(1)|^{2p}\right] \ar\le\ar  K\mbb{E}\left[\left|\int_0^{T\wedge \kappa_n}\<X_s, \gamma(\mu_s, \cdot)\> d s\right|^p \right]\cr
\ar\le\ar K + K\mbb{E}\left[\int_0^{T\wedge \kappa_n}\<X_s, 1\>^{2p} d s\right]\cr
\ar\le\ar K + K\int_0^{T}\mbb{E}\left[\sup_{0 \le t \le (s \wedge \kappa_n)}\<X_t, 1\>^{2p}\right] d s,
\eeqnn
where the constant $K$ depends on $T.$ Therefore, by \eqref{0413a} we have
\beqnn
\mbb{E}\left[\sup_{0 \le t \le T\wedge \kappa_n}\<X_t, 1\>^{2p}\right] \ar\le\ar K + \mbb{E} \left[\sup_{0 \le t \le T\wedge \kappa_n}|M_t(1)|^{2p}\right]\cr
\ar\le\ar K + K\int_0^{T}\mbb{E}\left[\sup_{0 \le t \le (s \wedge \kappa_n)}\<X_t, 1\>^{2p}\right] d s,
\eeqnn
which implies that $\mbb{E}\left[\sup_{0 \le t \le T\wedge \kappa_n}\<X_t, 1\>^{2p}\right] \le Ke^{KT}$
by Gronwall's inequality (see, e.g., \cite[Theorem 1]{WW64}). Letting $n \rightarrow \infty,$ the result follows from the monotone convergence theorem.
\qed
\begin{lemma}\label{Xp}
Let $(X_t)_{t \ge 0}$ satisfy the MP (\ref{0.1}, \ref{0.2}) and $T > 0$ be fixed. Then for any $t \in [0, T)$ and $p \ge 1$ we have
$\mbb{E}[\langle X_t, p_{T - t}(x - \cdot)\rangle] = \langle X_0, p_T(x - \cdot)\rangle$ and
\beqnn
\mbb{E}\left[\sup_{0 \le t \le T}\langle X_t, p_{T - t}(x - \cdot)\rangle^{2p}\right] \le \frac{K(T^{-p} + T^p)}{1 - 2^{-p}},
\eeqnn
where $\<X_T, p_0(x - \cdot)\> := \lim_{t \rightarrow T-} \<X_t, p_{T - t}(x - \cdot)\>.$
\end{lemma}
\proof
 By Lemma~\ref{t2.3} and the proof of \cite[Theorem~7.26]{L11}, we have
\beqlb\label{star}
\langle X_t, p_{T - t}(x - \cdot)\rangle = \langle X_0, p_T(x - \cdot)\rangle + \int_0^t\int_{\mbb{R}}p_{T - s}(x - z) M(d s, d z)
\eeqlb
for $t < T$, where $M(d s, d z)$ is the orthogonal martingale measure defined in Lemma~\ref{t2.3}. For any $t \in [0, T),$ by \eqref{Gamma1} and \eqref{p_t} we have
\beqlb\label{0404}
\ar\ar\mbb{E} \left[\int_0^t\int_\mbb{R} p_{T - s}(x - z)^2 \gamma(\mu_s, z)X_s(d z) d s\right] \cr
\ar\ar\qquad\le K \int_0^t \frac{ \mbb{E}\left[\langle X_s, 1\rangle\right]}{T - s}d s = K\<X_0, 1\> \ln \frac{T}{T - t} < \infty
\eeqlb
since $t \in [0, T)$. It implies that $t \mapsto \int_0^t \int_{\mbb{R}}p_{T - s}(x - z) M(d s, d z)$ is a martingale with $t \in [0,T).$ Then by \eqref{star} we have $\mbb{E}[\langle X_t, p_{T - t}(x - \cdot)\rangle] = \langle X_0, p_T(x - \cdot)\rangle$ for  $t < T$.
Further, by \eqref{p_t},
\beqlb\label{0409}
\langle X_0, p_T(x - \cdot)\rangle^{2p} \le T^{-p} \langle X_0, 1\rangle^{2p}
\eeqlb
 is bounded. Recall $\gamma$ is bounded. Let $0 < r_0 < T$ be fixed. Then, by Burkholder-Davis-Gundy's inequality and \eqref{p_t}, there are constants $K_1, K_2$ independent of $T, r_0$ such that
\beqlb\label{0409a}
I_{T,r_0} \ar:=\ar 2^{2p-1}\mbb{E}\left[\sup_{0 \le t \le T - r_0}\left|\int_0^t\int_{\mbb{R}} p_{T - s}(x - z)M(d s, d z)\right|^{2p}\right] \cr
\ar \le\ar  K_1 \mbb{E}\left[\left|\int_0^{T-r_0}\int_{\mbb{R}}p_{T - s}(x - z)^2 \gamma(\mu_s, z) X_s(d z)d s\right|^p\right] \cr
\ar \le\ar  K_2 \mbb{E}\left[\left|\int_0^{T-r_0}\int_{\mbb{R}}p_{T - s}(x - z)^2 X_s(d z)d s\right|^p\right]\cr
\ar\le\ar  \mbb{E}\left[\left|\int_0^{T-r_0}\frac{ K_2^{1/p}\langle X_s, p_{T - s}(x - \cdot)\rangle }{\sqrt{T - s}}d s\right|^p\right].
\eeqlb
Notice that $ab \le a^2 + b^2$ for $a, b \in \mbb{R}.$ We have
\beqnn
K_2^{1/p}\langle X_s, p_{T - s}(x - \cdot)\rangle \ar=\ar (2K_2^{1/p}T^{1/4})(2^{-1}T^{-1/4} \langle X_s, p_{T - s}(x - \cdot)\rangle)\\
\ar\le\ar 4K_2^{2/p}T^{1/2} + \frac{1}{4}T^{-1/2}\sup_{0 \le s \le T - r_0}\langle X_s, p_{T - s}(x - \cdot)\rangle^2.
\eeqnn
By the above inequality and \eqref{0409a} we have
\beqnn
I_{T,r_0} \le 8^p K_2^2 T^p + \frac{1}{2^p}\mbb{E}\left[\sup_{0 \le t \le T - r_0}\langle X_t, p_{T - t}(x - \cdot)\rangle^{2p}\right].
\eeqnn
By the above inequality, \eqref{star}, \eqref{0409} and using the fact $(a+b)^{2p} \le 2^{2p-1}a^{2p} + 2^{2p-1}b^{2p}$ for any $a, b \ge 0,$ we obtain
\beqlb\label{04051}
\ar\ar\mbb{E}\left[\sup_{0 \le t \le T-r_0}\langle X_t, p_{T - t}(x - \cdot)\rangle^{2p}\right]\cr
\ar\ar\quad \le 2^{2p-1} \langle X_0, p_T(x - \cdot)\rangle^{2p} + I_{T,r_0} \cr
\ar\ar\quad \le K(T^{-p} + T^p) + \frac{1}{2^p}\mbb{E}\left[\sup_{0 \le t \le T - r_0}\langle X_t, p_{T - t}(x - \cdot)\rangle^{2p}\right]
\eeqlb
with $K =(2^{2p-1}\langle X_0, 1\rangle^{2p})\vee (8^pK_2^2).$ By Lemma \ref{0405} and \eqref{p_t} one can check that
\beqnn
\mbb{E}\left[\sup_{0 \le t \le T-r_0}\langle X_t, p_{T - t}(x - \cdot)\rangle^{2p}\right] \le r_0^{-p} \mbb{E}\left[\sup_{0 \le t \le T-r_0}\langle X_t, 1\rangle^{2p}\right] < \infty.
\eeqnn
Moving the last term of \eqref{04051} to the left side of the inequality, we then have
\beqnn
\mbb{E}\left[\sup_{0 \le t \le T-r_0}\langle X_t, p_{T - t}(x - \cdot)\rangle^{2p}\right] \le \frac{K(T^{-p} + T^p)}{1 - 2^{-p}}.
\eeqnn
Letting $r_0 \rightarrow 0+$, by Fatou's lemma, we have
 \beqnn
\mbb{E}\left[\sup_{0 \le t \le T}\langle X_t, p_{T - t}(x - \cdot)\rangle^{2p}\right] \ar\le\ar \lim_{r_0\rightarrow 0+}\mbb{E}\left[\sup_{0 \le t \le T-r_0}\langle X_t, p_{T - t}(x - \cdot)\rangle^{2p}\right]\cr
\ar \le\ar \frac{K(T^{-p} + T^p)}{1 - 2^{-p}}.
\eeqnn
The result follows.
\qed

\begin{lemma}\label{joint_cont}
Let $M(d s, d z)$ be the martingale measure defined in Lemma~\ref{t2.3}. Then for all $t > 0 , x \in \mbb{R}$ the stochastic integral
\beqlb\label{Mtx}
M_t(x) = \int_0^{t}\int_{\mbb{R}}p_{t - s}(x - z)M(d s, d z)
\eeqlb
is well-defined. Moreover, $\{M_t(x): t > 0, x \in \mbb{R}\}$ has a modification which is almost surely continuous.
\end{lemma}
\proof
Recall that $\gamma$ is bounded by Condition~\ref{0.4b}. By \eqref{p_t} and Lemma \ref{Xp}, for every $t > 0$ and $x \in \mbb{R}$ we have
\beqnn
&\ &\mbb{E}\Big[\int_0^{t} \int_{\mbb{R}}p_{t - s}(x - z)^2 {\gamma(\mu_s, z)}X_s(d z)d s\Big]\\
&\ &\qquad \le  K \mbb{E}\Big[\int_0^{t} \int_{\mbb{R}}p_{t - s}(x - z)^2 X_s(d z)d s\Big]\\
&\ &\qquad \le K\int_0^{t} \frac{\mbb{E}\left[\<X_s, p_{t - s}(x - \cdot)\>\right]}{\sqrt{t - s}} d s \\
&\ &\qquad \le K \langle X_0, p_t(x - \cdot)\rangle\sqrt{t} \le K\langle X_0, 1\rangle,
\eeqnn
which is bounded. Then $M_t(x)$ is well-defined for all $t, x$ (see, e.g., \cite[p. 55]{IW89}). Let $0 < r_0 < T$ be fixed. For any $0 < r_0 \le r \le t \le T$ and $p \ge 1,$ we have
\beqlb\label{I_1I_2}
\ar\ar\mbb{E}[|M_t(x) - M_r(x)|^{2p}]\cr
\ar \ar\quad\le K\mbb{E}\left[\left| \int_r^t \int_{\mbb{R}}p_{t - s}(x - z)M(d s, d z)\right|^{2p}\right]\cr
\ar\ar\qquad  + K\mbb{E}\left[\left| \int_0^{r} \int_{\mbb{R}}[p_{r - s}(x - z) - p_{t - s}(x - z)]M(d s, d z)\right|^{2p}\right]\cr
\ar\ar\quad=: I_1 + I_2.
\eeqlb
 Notice that $a \le a^2 +1$ for $a \in \mbb{R}.$  By Lemma~\ref{Xp}, for any $r_0 \le t \le T$ we have
\beqlb\label{add1}
\mbb{E}\left[\sup_{0 \le s \le t}\langle X_s, p_{t - s}(x - \cdot)\rangle^{p}\right] \ar\le\ar \mbb{E}\left[\sup_{0 \le s \le t}\langle X_s, p_{t - s}(x - \cdot)\rangle^{2p}\right]  +1\cr
\ar\le\ar \frac{K(t^{-p} + t^p)}{1 - 2^{-p}} + 1\cr
\ar\le\ar \frac{K(r_0^{-p} + T^p)}{1 - 2^{-p}} + 1,
\eeqlb
where the above $K$ is independent of $t.$ By \eqref{p_t}, \eqref{add1} and Burkholder-Davis-Gundy's inequality one can see that
\beqlb\label{I_11}
I_1 \ar\le\ar K\mbb{E}\left[\left| \int_r^t  d s \int_{\mbb{R}}p_{t - s}(x - z)^2X_s(d z)\right|^{p}\right]\cr
\ar\le\ar K\mbb{E}\left[\left| \int_r^t  \frac{\langle X_s, p_{t - s}(x - \cdot)\rangle}{\sqrt{t - s}}d s \right|^{p}\right]\cr
\ar\le\ar K(t - r)^{p/2}\mbb{E}\left[\sup_{0 \le s \le t}\langle X_s, p_{t - s}(x - \cdot)\rangle^p\right] \le K(t - r)^{p/2},
\eeqlb
where the above constant $K$ only depends on $r_0, T.$ In the following of the proof we assume that $\delta \in (0, 1/4).$ By \cite[Lemma~III 4.5]{P02}, for any $0 < s < r < t \le T$ we have
\beqlb\label{p_{r - s}}
\ar\ar\left[p_{r - s}(x - z) - p_{t - s}(x - z)\right]^2\cr
\ar\ar\quad \le (t - r)^{\delta}(r - s)^{-3\delta/2}\left[p_{r - s}(x - z)^{2 - \delta} + p_{t - s}(x - z)^{2 - \delta}\right].
\eeqlb
It follows from \eqref{p_{r - s}} and Burkholder-Davis-Gundy's inequality that
\beqnn
 I_2 \ar\le\ar K\mbb{E}\left[\left| \int_0^{r} d s\int_{\mbb{R}}[p_{r - s}(x - z) - p_{t - s}(x - z)]^2X_s(d z)\right|^{p}\right]\cr
\ar\le\ar K (t\! -\! r)^{p\delta}\mbb{E}\left[\left| \int_0^r\! (r \!-\! s)^{-3\delta/2}\< X_s, p_{r - s}(x\! -\! \cdot)^{2 - \delta} + p_{t - s}(x\! - \!\cdot)^{2 - \delta}\> d s\right|^{p}\right].
\eeqnn
 Moreover, by \eqref{p_t} and \eqref{add1} once again we have
\beqnn
\ar\ar\mbb{E}\left[\left| \int_0^{r} (r - s)^{-3\delta/2}\left\langle X_s, p_{r - s}(x - \cdot)^{2 - \delta} + p_{t - s}(x - \cdot)^{2 - \delta}\right\rangle d s\right|^{p}\right]\cr
\ar\ar\quad \le K\mbb{E}\left[\sup_{0 \le s \le r}\langle X_s, p_{r - s}(x - \cdot)\rangle^p\right] \left|\int_0^r (r - s)^{-(1+2\delta)/2}d s\right|^p\cr
\ar\ar\qquad + K \mbb{E}\left[\sup_{0 \le s \le t}\langle X_s, p_{t - s}(x - \cdot)\rangle^p\right] \left|\int_0^r (r - s)^{-3\delta/2} (t - s)^{-(1-\delta)/2}ds\right|^p\\
\ar\ar\quad \le K\mbb{E}\left[\sup_{0 \le s \le r}\langle X_s, p_{r - s}(x - \cdot)\rangle^p\right] r^{(1 - 2\delta)p/2} \cr
\ar\ar\qquad +K \mbb{E}\left[\sup_{0 \le s \le t}\langle X_s, p_{t - s}(x - \cdot)\rangle^p\right]\left(\int_0^r (r - s)^{-3\delta} ds\int_0^r(t - s)^{\delta - 1}d s \right)^{p/2} \cr
\ar\ar\quad \le K\left(\frac{K(r_0^{-p} + T^p)}{1 - 2^{-p}} + 1\right) T^{(1 - 2\delta)p/2}.
\eeqnn
By the above two inequalities, we have
\beqlb\label{I_22}
I_2 \le K (t - r)^{p\delta}, \qquad 0 < r_0 \le r \le t \le T
\eeqlb
with the above constant $K$ only depending on $r_0, T.$ Combining \eqref{I_1I_2}, \eqref{I_11} and \eqref{I_22} we obtain
\beqlb\label{M_tx}
\mbb{E}[|M_t(x) - M_r(x)|^{2p}] \le K [(t - r)^{p/2} + (t - r)^{p\delta}]
\eeqlb
for any $0 < r_0 \le r \le t \le T$, where the above constant $K$ only depends on $r_0, T.$

On the other hand, by \cite[(2.4e)]{R87}, for any $0 < \beta < 1$ we have
\beqnn
|p_{t - s}(x - z) - p_{t - s}(y - z)| \le K |x - y|^\beta (t - s)^{-(1+\beta)/2}, \quad  x, y \in \mbb{R}.
\eeqnn
Then for all $x, y \in \mbb{R},$ by Burkholder-Davis-Gundy's inequality,  \eqref{add1} and Lemma \ref{Xp} one can check that, for $0 < r_0 \le t \le T,$
\beqlb\label{Mtxy}
\ar\ar\mbb{E}[|M_t(x) - M_t(y)|^{2p}] \cr
\ar\ar\quad= \mbb{E}\left[\left|\int_0^{t}\int_{\mbb{R}}[p_{t - s}(x - z) - p_{t - s}(y - z)]M(d s, d z)\right|^{2p}\right]\cr
\ar\ar\quad\le K\mbb{E}\left[\left(\int_0^{t}\int_{\mbb{R}}[p_{t - s}(x - z) - p_{t - s}(y - z)]^2X_s(d z)d s\right)^{p}\right]\cr
 \ar\ar\quad \le K|x - y|^{\beta p}\mbb{E}\!\left[\!\left(\!\int_0^{t} \!\!(t - s)^{-(1+\beta)/2}\!\<X_s, p_{t - s}(x - \cdot) + p_{t - s}(y - \cdot)\>\! d s\right)^{p}\!\right]\cr
\ar\ar\quad\le K|x - y|^{\beta p}\left(\int_0^t (t - s)^{-(1+\beta)/2}d s\right)^p\cr
\ar\ar\qquad\cdot\left\{\mbb{E}\left[\sup_{0 \le s \le t}\langle X_s, p_{t - s}(x - \cdot)\rangle^p\right] + \mbb{E}\left[\sup_{0 \le s \le t}\langle X_s, p_{t - s}(y - \cdot)\rangle^p\right]\right\}\cr
\ar\ar\quad \le K|x - y|^{\beta p}
\eeqlb
with $\beta \in (0, 1),$ where the constant $K$ depends on $r_0, T.$ By \eqref{M_tx} and \eqref{Mtxy},
\beqnn
\mbb{E} [|M_t(x) - M_r(y)|^{2p}] \le K[(t - r)^{p/2} + (t - r)^{p\delta} + |x - y|^{\beta p}],
\eeqnn
for $0 < r_0 \le r \le t \le T$ and $x, y \in \mbb{R},$ where $p$ can be taken greater than $\max(1/\delta,1/\beta)$. By Kolmogorov's continuity criteria (see, e.g., \cite[Theorem~3.23]{K02}), $\{M_t(x): t \in [r_0, T], x \in \mbb{R}\}$ has a continuous version. The result follows by taking suitable sequences $\{r_0, r_1, \cdots\}$ and $\{T_1, T_2, \cdots\}$ such that $\lim_{n \rightarrow \infty}r_n = 0$ and $\lim_{n \rightarrow \infty}T_n = +\infty.$
\qed

\begin{prop}\label{t2.4}
Suppose that $(X_t)_{t \ge 0}$ is a solution to the MP (\ref{0.1}, \ref{0.2}). Then almost surely for each $t \ge 0,$ the random measure $X_t(d x)$ is absolutely continuous respect to $d x$ with density $\mu_t(x)$ satisfying \eqref{1.1}.
Conversely, assume $(\mu_t)_{t \ge 0}$ is a solution to \eqref{1.1}. Then there exists a solution to the MP (\ref{0.1}, \ref{0.2}).
\end{prop}

\proof
Suppose that $(X_t)_{t \ge 0}$ satisfies the MP (\ref{0.1}, \ref{0.2}). By Lemma~\ref{t2.3} and the proof of \cite[Theorem~7.26]{L11}, for any $\phi \in B(\mbb{R})$ we have
\beqnn
\<X_t, \phi\> = \<X_0, P_t\phi\> + \int_0^t\int_{\mbb{R}}P_{t - s}\phi(x)M(d s, d x),
\eeqnn
where $M(d s, d x)$ is the orthogonal martingale measure defined in Lemma~\ref{t2.3} and $P_t$ is the semigroup generated by $\Delta/2$. Taking an appropriate modification, we may and shall assume that $\{M_t(x): t > 0, x \in \mbb{R}\}$ is almost surely continuous by Lemma \ref{joint_cont}. Then by stochastic Fubini's theorem (e.g., see Li~\cite[Theorem~7.24]{L11}), for every $t \in (0, T]$ and $\phi \in C_c(\mbb{R})$ we get $\<X_t, \phi\> =  \int_{\mbb{R}} \mu_t(x)\phi(x)d x$ almost surely with
\beqnn
\mu_t(x) = \int_\mbb{R}p_t(x - z)\mu_0(z)d z + M_t(x).
\eeqnn
 Note that $t \mapsto \<X_t, \phi\>$ and $t \mapsto
 \int_{\mbb{R}} \mu_t(x)\phi(x)d x$ are continuous. Considering suitable sequences $\{\phi_1, \phi_2, \cdots\} \subset C_c(\mbb{R})$ and $\{t_1, t_2, \cdots\} \subset (0, \infty)$, by Lemma~\ref{joint_cont} we get almost surely for all $t > 0$ and $\phi \in C_b^2(\mbb{R})$, $\<X_t, \phi\> = \int_{\mbb{R}}\mu_t(x)\phi(x)d x,$ i.e., $X_t(d x)$ is absolutely continuous respect to Lebesgue measure $d x.$ Set $p_t(x - z) = 0$ for all $t < 0.$ By El Karoui and M\'{e}l\'{e}ard \cite[Theorem III-6]{EM90}, on some extension of the probability space one can define a white noise $W(d s, d z)$ on $\mbb{R}_+\times \mbb{R}$ based on $d s d z$ such that almost surely for any $u > 0$, $t > 0$ and $x \in \mbb{R},$
\beqnn
\int_0^u\int_{\mbb{R}}p_{t - s}(x - z)M(d s, d z) = \int_0^u\int_{\mbb{R}}\sqrt{\mu_s(z)\gamma(\mu_s, z)}p_{t - s}(x - z)W(d s, d z),
\eeqnn
which implies that
\beqlb\label{mu}
\mu_t(x) = \<\mu_0, p_t(x - \cdot)\>  + \int_0^t\int_{\mbb{R}}\sqrt{\mu_s(z)\gamma(\mu_s, z)}p_{t - s}(x - z)W(d s, d z)
\eeqlb
for all $t > 0, x \in \mbb{R}$ almost surely, and $\mu_t(x)$ is joint continuous in $(t, x) \in (0, \infty)\times \mbb{R}.$ Notice that
 \beqnn
 \int_{\mbb{R}}\langle f, p_t(x - \cdot)\rangle^2d x \le K\int_{\mbb{R}}d x\int_\mbb{R}f(y)^2p_t(x - y)d y \le K\int_\mbb{R}f(y)^2d y < \infty
 \eeqnn
 for any $f \in C_c(\mbb{R})^+.$ Then we have $\lim_{t \rightarrow 0}\int_{\mbb{R}}[\<f, p_t(x - \cdot)\> -f(x)]^2d x = 0$ by dominated convergence theorem. Moreover, by Lemma \ref{Xp} we have
\beqnn
\ar\ar\int_{\mbb{R}} d x\mbb{E}\left[\int_0^t\int_{\mbb{R}}\mu_s(z)p_{t - s}(x - z)^2 d s d z\right]\cr
\ar\ar\quad \le K \left[\int_{\mbb{R}} d x\int_0^t\frac{\mbb{E}[\<\mu_s,p_{t - s}(x - \cdot)\>]}{\sqrt{t - s}} d s\right]\cr
\ar\ar\quad  \le K\sqrt{t} \int_{\mbb{R}} \<\mu_0, p_t(x - \cdot)\>d x \rightarrow 0, \qquad t \rightarrow 0.
\eeqnn
It then implies that $\mu_t \rightarrow \mu_0$ in $L^2(\mbb{R})$ as $t \rightarrow 0$ in the following sense:
	\beqlb\label{0730}
	\mbb{E}\left[\|\mu_t - \mu_0\|^2_{L^2}\right]\! \!\ar=\ar\!\! \mbb{E}\left[\int_{\mbb{R}}|\mu_t(x) - \mu_0(x)|^2 d x\right]\cr
	\ar\le\ar\!\! K \int_\mbb{R}|\langle \mu_0, p_t(x - \cdot)\rangle - \mu_0(x)|^2 d x\cr
	\ar\ar\!  + K \int_\mbb{R}\!d x\mbb{E}\left[\int_0^t\!\int_{\mbb{R}}\mu_s(x)p_{t - s}(x - z)^2 d z d s \right]\cr
	\ar\rightarrow\ar\! 0
	\eeqlb
as $ t \rightarrow 0.$ Then \eqref{mu} holds for all $t \ge 0, x \in \mbb{R}$ almost surely. 
And \eqref{1.1} holds similar to \cite[Theorem~2.1]{S94}.

Conversely, suppose that $(\mu_t)_{t \ge 0}$ satisfies \eqref{1.1} and  we denote $X_t(d x) = \mu_t(x)d x$. It then follows from Lemma~\ref{<mu,1>} that $X_t\in M(\mbb{R})$ almost surely for every $t \ge 0.$ For any $\phi \in C_b^2(\mbb{R})$ one can check that
\beqnn
\<X_t, \phi\> \ar=\ar \int_{\mbb{R}}\mu_t(x)\phi(x) d x\cr
 \ar=\ar \int_{\mbb{R}}\mu_0(x)\phi(x) d x + \frac{1}{2}\int_0^t \int_{\mbb{R}}\mu_s(x)\phi''(x) d x d s\cr
\ar\ar +\int_0^t\int_{\mbb{R}}\phi(x)\sqrt{\mu_s(x)\gamma(\mu_s, x)}W(d s, d x)\cr
\ar=\ar \<X_0, \phi\> + \frac12\int_0^t\<X_s, \phi''\>d s + M_t(\phi),
\eeqnn
 where $(M_t(\phi))_{t \ge 0}$ is a continuous local martingale with quadratic variation process $(\<M(\phi)\>_t)_{t \ge 0}$ satisfying \eqref{0.2}. Further, by Lemma \ref{<mu,1>}, for $\phi \in C_b^2(\mbb{R}),$
\beqnn
\mbb{E}[\<M(\phi)\>_t] \le K \int_0^t \mbb{E} [\langle\mu_s, 1\rangle] d s  = K t\langle\mu_0, 1\rangle,
\eeqnn
 which implies $(M_t(\phi))_{t \ge 0}$ is a martingale (see \cite[p. 55]{IW89}). The proof ends.
\qed


Now we show the existence of the solution to \eqref{1.1}. For any $T > 0,$ let $m \ge 1$ and $t_k = kT/m$ with $k = 0, 1, \cdots, m.$ For $x \in \mbb{R},$ let $\rho$ be the mollifier given by
\beqlb\label{rho}
\rho(x) := C\exp\{-1/(1 - x^2)\}1_{\{|x| < 1\}},
\eeqlb
 and $C$ is a constant such that $\int_{\mbb{R}}\rho(x) d x = 1.$
 Let $\rho_m = m\rho(mx)$ and $\gamma_m(\cdot, x) = \int_{\mbb{R}} \rho_{m}(x - y)\gamma(\cdot, y)d y$  be the mollification of $\gamma$. Then $x \mapsto \gamma_m(\cdot, x)$ is continuous and $\lim_{m \rightarrow \infty}\gamma_m(\cdot, x) = \gamma(\cdot, x)$ almost every $x\in\mbb{R}$ (see, e.g., \cite[p. 630, Theorem 6]{E98}). Further, by Condition \ref{0.4b} one can see that
 \beqlb\label{gamma0722}
 \sup_{x \in \mbb{R}, m \ge 1, \atop f \in L^1(\mbb{R})^+}\gamma_m(f, x) \ar=\ar \sup_{x \in \mbb{R}, m \ge 1, \atop f \in L^1(\mbb{R})^+}\int_{\mbb{R}} \rho_{m}(x - y)\gamma(f, y)d y\cr \ar\le\ar\sup_{x \in \mbb{R}, m \ge 1}  K \int_{\mbb{R}} \rho_{m}(x - y)d y = K.
 \eeqlb 
For any $f \in C_c^\infty(\mbb{R}),$ we define a sequence of approximation by
\beqlb\label{mun}
\<\mu_t^{m}, f\> \ar=\ar \<\mu_0, f\> + \frac12\int_0^t\<\mu_s^m, f''\> d s \cr
\ar\ar\!+\! \sum_{k = 1}^m\int_{ t_{k - 1}\wedge t}^{ t_k\wedge t}\! \int_{\mbb{R}}\! \sqrt{\gamma_m(\mu_{t_{k\! -\! 1}}^m\!, x)} G_m(\mu_s^m(x))f(x)W(d s, d x),
\eeqlb
where
\beqnn
G_m(x) = \int_{\mbb{R}} \left[p_{m^{-1}}(x - y) - p_{m^{-1}}(y)\right](\sqrt{|y|}\wedge m)d y
\eeqnn
is a Lipschitz function for fixed $m \ge 1.$ Moreover, one can see that $G_m(0) = 0$, $\lim_{m \rightarrow \infty}G_m(x) = \sqrt{x}$ for all $x \ge 0$. 

Recall that $C([0, T], C_{\rm tem}^+(\mbb{R}))$ is the collection of all continuous maps from $[0, T]$ to $C_{\rm tem}^+(\mbb{R})$ with topology induced by \eqref{d1fg} by changing $\mathcal{Y}$ to $C_{\rm tem}^+(\mbb{R})$. And $C([0, T]\times \mbb{R}, \mbb{R}_+)$ is the space of continuous maps from $[0,T]\times \mbb{R}$ to $\mbb{R}_+$ furnished with metric \eqref{dfg}. The topology of $C([0, T], C_{\rm tem}^+(\mbb{R}))$ is stronger than that of $C([0, T]\times \mbb{R}, \mbb{R}_+)$. Hence $C([0, T], C_{\rm tem}^+(\mbb{R}))$ is a subspace of $C([0, T]\times \mbb{R}, \mbb{R}_+)$. Then we have the following result.
 \begin{lemma}	
 For every $m \ge 1,$ there is a pathwise unique continuous $C_{\rm tem}^+(\mbb{R})$-valued solution $(\mu_t^m)_{t \in [0, T]}$ to \eqref{mun}. Furthermore,  $\{\mu_t^{m}(x): t \in [0, T], x \in \mbb{R}\} \in C([0, T] \times \mbb{R}, \mbb{R}_+)$ satisfying, for every $\lambda > 0$,
	\beqlb\label{0719}
\sup_{0 \le t \le T}\sup_{x \in \mbb{R}}\left[e^{-\lambda |x|}\mu_t^m(x)\right] < \infty \qquad a.s.
		\eeqlb 
\end{lemma}
\proof
For every $m \ge 1,$ there is a constant $K > 0$ such that
\beqlb\label{G_m}
|G_m(x)| \ar\le\ar K + \int_{\mbb{R}}p_{m^{-1}}(x - y)\sqrt{|y|}d y = K + \mbb{E}\left[\sqrt{ |W_{m^{-1}}^x|}\right]\cr
\ar \le\ar K +  \left[\mbb{E}\left[ |W_{m^{-1}}^x|^2\right]\right]^{1/4} \le K + (1 + x^2)^{1/4}\cr
\ar \le\ar K(\sqrt{|x|} + 1),
\eeqlb
where $(W_t^x)_{t \ge 0}$ is a Brownian motion with initial value $x.$
 Further, $|G_m(x)| \le K(|x| +1).$  Recall $\mu_0 \in C_c(\mbb{R})^+$, $x \mapsto \gamma_m(\cdot, x)$ is continuous and bounded. And $G_m(x)$ is Lipschitz with $G_m(0) = 0.$ By \cite[Theorems~2.2 and 2.3]{S94}, when $k = 1,$ there is a pathwise unique continuous $C_{\rm tem}^+(\mbb{R})$-valued solution $(\mu_t^{m})_{t \in [t_{k - 1}, t_k]}$ to
 \beqlb\label{11221}
 \<\mu_t^{m}, f\> \ar=\ar \<\mu_{t_{k - 1}}^m, f\> + \frac12\int_{t_{k - 1}}^t\< \mu_s^m, f''\> d s\cr
 \ar\ar + \int_{t_{k - 1}}^t\int_{\mbb{R}} \sqrt{\gamma_m(\mu_{t_{k - 1}}^m, x)} G_m(\mu_s^m(x))f(x)
 W(d s, d x), 
 \eeqlb
for $t \in [t_{k - 1}, t_k]$ and any $f \in C_c^\infty(\mbb{R})$ almost surely, which implies that
 \beqlb\label{11222}
 \sup_{t \in [t_{k - 1}, t_k]}\sup_{x \in \mbb{R}}[e^{-\lambda |x|}\mu_t^m(x)] < \infty
 \eeqlb
 almost surely for each $\lambda > 0$.  Notice that $C([t_{k - 1}, t_k], C_{\rm tem}^+(\mbb{R}))$ is a subspace of $C([t_{k - 1}, t_k] \times \mbb{R}, \mbb{R}_+)$, which implies that  $\{\mu_t^{m}(x): t \in [t_{k - 1}, t_k], x \in \mbb{R}\} \in C([t_{k - 1}, t_k] \times \mbb{R}, \mbb{R}_+)$. Conditioned on $\mcr{F}_{t_{k - 1}},$ the above result still holds for $k = 2, \cdots, m$ by mathematical induction. In other word, for each $k = 1, \cdots, m,$ there is a pathwise unique continuous $C_{\rm tem}^+(\mbb{R})$-valued $(\mu_t^{m})_{t \in [t_{k - 1}, t_k]}$ to \eqref{11221}. Moreover,  $\{\mu_t^{m}(x): t \in [t_{k - 1}, t_k], x \in \mbb{R}\} \in C([t_{k - 1}, t_k] \times \mbb{R}, \mbb{R}_+)$ satisfying \eqref{11222}. The result follows.
\qed

For $x \in \mbb{R},$ let $J(x) = \int_{\mbb{R}}e^{-|y|}\rho(x - y) d y,$
where $\rho$ is the mollifier given by \eqref{rho}. Let $J^{(n)}(x)$ be the $n$th derivative of $J(x).$ By Mitoma~\cite[(2.1)]{M85}, for $n \ge 0,$ there are constants $c_n, C_n$ such that
\beqlb\label{J}
c_n e^{-|x|} \le J^{(n)}(x) \le C_n e^{-|x|}, \qquad  x \in \mbb{R}.
\eeqlb

\begin{lemma}\label{l2.4}
 For every $T > 0$ and $p \ge 1$, we have
\beqnn
\sup_{0 \le t \le T, m \ge 1}\mbb{E}\left[\int_{\mbb{R}}\mu_t^m(x)^{2p}J(x)d x\right] < \infty.
\eeqnn
\end{lemma}
\proof
Let $0 \le t \le T.$ Using the convolution form, the solution $(\mu_t^m)_{t \ge 0}$ to \eqref{mun} can be represented as
{\small\beqlb\label{mun1}
\mu_t^m(x) \ar=\ar \<\mu_0, p_t(x - \cdot)\>\ccr
\ar\ar + \sum_{k = 1}^m\!\int_{t_{k - 1}\wedge t}^{t_k\wedge t}\!\int_{\mbb{R}}\! \sqrt{\gamma_m(\mu_{t_{k - 1}}^m\!, z) }G_m(\mu_s^m(z))p_{t - s}(x\! -\! z) W(d s, d z).
\eeqlb}
 For each $n \ge 1$ we define stopping time
$\sigma_n^m = \inf\left\{t \ge 0:\! \int_{\mbb{R}}\mu_t^m(x)^{2p}J(x)d x\! \ge\! n\right\}.$ Recall that $\mu_0 \in C_c(\mbb{R})^+.$ For any $p \ge 1$, it follows from \eqref{0719} that
\beqnn
\ar\ar\sup_{0 \le t \le T}\int_{\mbb{R}}\mu_t^m(x)^{2p}J(x)d x\cr
\ar\ar\qquad\le \left[\sup_{0 \le t \le T}\sup_{x \in \mbb{R}}\left[e^{- |x|/4p}\mu_t^m(x)\right]\right]^{2p}\int_{\mbb{R}}e^{|x|/2}J(x)d x < \infty
\eeqnn
almost surely for every $T > 0$ and $m \ge 1$. We then have $\sigma_n^m \rightarrow \infty$ almost surely as $n \rightarrow \infty$. By \eqref{p_t} and H{\"o}lder's inequality, for any $t \le T$ we have
\beqlb\label{hh}
\ar\ar\left|\int_0^t 1_{\{s \le \sigma_n^m\}}d s\int_{\mbb{R}} \mu_s^m(z)p_{t - s}(x - z)^2 d z\right|^{2p}\cr
\ar\ar\qquad\le \left|\int_0^t \frac{1_{\{s \le \sigma_n^m\}}}{\sqrt{t - s}}d s\int_{\mbb{R}} \mu_s^m(z)p_{t - s}(x - z) d z\right|^{2p}\cr
\ar\ar \qquad\le K\left[\int_0^t \frac{1_{\{s \le \sigma_n^m\}}}{\sqrt{t - s}}\left|\int_{\mbb{R}} \mu_s^m(z)p_{t - s}(x - z) d z\right|^{2p}d s\right] \cdot \left[\left|\int_0^t \frac{1}{\sqrt{t - s}}d s\right|^{2p/q}\right]\cr
\ar\ar \qquad\le K\int_0^t \frac{1_{\{s \le \sigma_n^m\}}}{\sqrt{t - s}}d s\int_{\mbb{R}} \mu_s^m(z)^{2p}p_{t - s}(x - z) d z
\eeqlb
with $1/2p + 1/q = 1$ and $p, q \ge 1,$ and the above $K$ depending on $T.$  By the definition of $J(\cdot)$, we have $J(u + z) \le J(z)e^{|u|}$ for $u, z \in \mbb{R}$. Moreover,
\beqlb\label{hhhh}
\int_{\mbb{R}}e^{|u|}p_{t - s}(u)d u \le \int_{\mbb{R}}e^{\sqrt{T}|u|}p_1(u)d u =  \mbb{E}\left[e^{\sqrt{T}|W_1|}\right] \le K
\eeqlb
for $0 \le s \le t \le T,$ where $(W_t)_{t \ge 0}$ is a standard Brownian motion. By \eqref{hh}, \eqref{hhhh} and a change of variable, we obtain
\beqlb\label{5h}
\ar\ar\mbb{E}\left[\int_{\mbb{R}}J(x)d x\left|\int_0^t 1_{\{s \le \sigma_n^m\}}d s \int_{\mbb{R}} \mu_s^m(z)p_{t - s}(x - z)^2 d z\right|^{2p}\right]\cr
\ar\ar\qquad \le K \mbb{E}\left[\int_0^t \frac{1_{\{s \le \sigma_n^m\}}}{\sqrt{t - s}} d s \int_{\mbb{R}}J(x)d x\int_{\mbb{R}}\mu_s^m(z)^{2p}p_{t - s}(x - z) d z\right] \cr
\ar\ar\qquad \le K \mbb{E}\left[\int_0^t \frac{1_{\{s \le \sigma_n^m\}}}{\sqrt{t - s}} d s \int_{\mbb{R}}\mu_s^m(z)^{2p} d z\int_{\mbb{R}}J(u + z) p_{t - s}(u) d u\right]\cr
\ar\ar\qquad \le K \mbb{E}\left[\int_0^t \frac{1_{\{s \le \sigma_n^m\}}}{\sqrt{t - s}} d s \int_{\mbb{R}}\mu_s^m(z)^{2p} J(z)d z\int_{\mbb{R}}e^{|u|} p_{t - s}(u) d u\right]\cr
\ar\ar\qquad \le K \int_0^t \frac{1}{\sqrt{t - s}}\mbb{E}\left[\int_{\mbb{R}}1_{\{s \le \sigma_n^m\}}\mu_s^m(z)^{2p} J(z) d z\right] d s,
\eeqlb
where the above $K$ depends on $T$ only. By \eqref{p_t} again we have
\beqlb\label{hhh}
\ar\ar\mbb{E}\left[\int_{\mbb{R}}J(x)d x\left|\int_0^t 1_{\{s \le \sigma_n^m\}}d s \int_{\mbb{R}} p_{t - s}(x - z)^2 d z\right|^p\right]\cr
\ar\ar\qquad \le \int_{\mbb{R}}J(x)d x\left|\int_0^t \frac{1}{\sqrt{t - s}}d s\right|^p  \le K,
\eeqlb
where again the constant $K$ only depends on $T.$ Now we denote
\beqnn
M_n(t,\! x)\! = \!\sum_{k = 1}^m\!\int_{t_{k - 1}\wedge t}^{t_k\wedge t}\!\int_{\mbb{R}}\! 1_{\{s \le \sigma_n^m\}}  \sqrt{\gamma_m(\mu_{t_{k - 1}}^m\!, z) } G_m(\mu_s^m(z))p_{t - s}(x\! -\! z) W(d s, d z).
\eeqnn
By Burkholder-Davis-Gundy's inequality, \eqref{gamma0722}, \eqref{G_m}, \eqref{5h} and \eqref{hhh}, we have
\beqlb\label{um1}
\ar\ar\int_{\mbb{R}}J(x)\mbb{E}\left[|M_n(t, x)|^{4p}\right]d x\cr
\ar\ar\qquad=\int_{\mbb{R}}J(x)d x\mbb{E}\Bigg[\Bigg|\int_0^t\int_{\mbb{R}} \sum_{k = 1}^m1_{(t_{k - 1}, t_k]}(s)1_{\{s < \sigma_n^m\}} \sqrt{\gamma_m(\mu_{t_{k - 1}}^m, z)} \cr
\ar\ar\qquad\qquad\qquad\qquad\qquad\qquad \cdot G_m(\mu_s^m(z))p_{t - s}(x - z) W(d s, d z)\Bigg|^{4p}\Bigg]\cr
\ar\ar\qquad \le K \int_{\mbb{R}}J(x)d x\mbb{E}\left[\left|\int_0^t \int_{\mbb{R}} 1_{\{s \le \sigma_n^m\}}G_m(\mu_s^m(z))^2p_{t - s}(x - z)^2 d s d z\right|^{2p}\right]\cr
\ar\ar\qquad \le K \mbb{E}\left[\int_{\mbb{R}}J(x)d x\left|\int_0^t d s \int_{\mbb{R}}1_{\{s \le \sigma_n^m\}} (\mu_s^m(z) + 1)p_{t - s}(x - z)^2 d z\right|^{2p}\right]\cr
\ar\ar\qquad \le K +  K \int_0^t \frac{1}{\sqrt{t - s}}\mbb{E}\left[\int_{\mbb{R}}1_{\{s \le \sigma_n^m\}}\mu_s^m(z)^{2p} J(z) d z\right] d s,
\eeqlb
 where the above constants $K$ only depends on $T.$ Recall that $\mu_0 \in C_c(\mbb{R})^+,$ which implies $\mu_0$ is bounded. Then we have
\beqlb\label{um2}
\int_{\mbb{R}}\<\mu_0, p_t(x - \cdot)\>^{2p}J(x)d x \ar=\ar \int_{\mbb{R}}J(x)d x\left[\int_{\mbb{R}}p_t(x - z)\mu_0(z)d z\right]^{2p}\cr
\ar\le\ar K\int_{\mbb{R}}J(x)d x,
\eeqlb
which is bounded. Note that
\beqlb\label{0409b}
\mu_t^m(x)1_{\{t \le \sigma_n^m\}} = \<\mu_0, p_t(x - \cdot)\>1_{\{t \le \sigma_n^m\}} + 1_{\{t \le \sigma_n^m\}}M_{n}(t, x).
\eeqlb
By H\"{o}lder's inequality and $a \le a^2 + 1$ for every $a \in \mbb{R},$ we then have
\beqnn
\mbb{E}\left[|1_{\{t \le \sigma_n^m\}}M_n(t, x)|^{2p} \right] \le \left\{\mbb{E}\left[|M_n(t, x)|^{4p}\right]\right\}^{1/2}
\le\mbb{E}\left[|M_n(t, x)|^{4p}\right] + 1.
\eeqnn
Combining the above inequality, \eqref{um1}, \eqref{um2} and \eqref{0409b}, we obtain
\beqnn
\ar\ar\mbb{E}\left[\int_{\mbb{R}}1_{\{t \le \sigma_n^m\}}\mu_t^m(x)^{2p}J(x)d x\right]\cr
\ar\ar\quad \le K + K\int_0^t \frac{1}{\sqrt{t - s}} \mbb{E}\left[\int_{\mbb{R}}1_{\{s \le \sigma_n^m\}}\mu_s^m(x)^{2p}J(x)d x\right]d s.
\eeqnn
Iterating the above once, one can check that
{\small\beqnn
\ar\ar\mbb{E}\left[\int_{\mbb{R}}1_{\{t \le \sigma_n^m\}}\mu_{t}^m(x)^{2p}J(x)d x\right]\cr
\ar \ar\quad\le K + K \mbb{E}\left[\int_0^{t}d r\int_{\mbb{R}}1_{\{r \le \sigma_n^m\}}\mu_r^m(x)^{2p}J(x)d x\int_r^{t} \frac{1}{\sqrt{({t} - s)(s - r)}}d s\right]\\
\ar\ar\quad\le K + K\int_0^t \mbb{E}\left[\int_{\mbb{R}}1_{\{r \le \sigma_n^m\}}\mu_{r}^m(x)^{2p}J(x)d x\right]d r.
\eeqnn}
By Gronwall's inequality (see, e.g., \cite[Theorem 1]{WW64}), we have
\beqnn
\mbb{E}\left[\int_{\mbb{R}}1_{\{t \le \sigma_n^m\}}\mu_{t}^m(x)^{2p}J(x)d x\right] \le Ke^{Kt}, \qquad 0 \le t \le T,
\eeqnn
where the above $K$ is independent of $t, n, m.$ Letting $n \rightarrow \infty,$ the result follows from the monotone convergence theorem.
\qed

We proceed to proving the tightness of $(\mu^m_\cdot)$ in $C([0, T] \times \mbb{R}, \mbb{R}_+).$ Denote
\beqlb\label{num}
\nu_t^m(x) = \sum_{k = 1}^m\int_{t_{k - 1}\wedge t}^{t_k\wedge t}\int_{\mbb{R}} \sqrt{\gamma_m(\mu_{t_{k - 1}}^m, z)} G_m(\mu_s^m(z))p_{t - s}(x - z) W(d s, d z).
\eeqlb

 \begin{lemma}\label{l22.6}
 For any fixed $0 < \delta < 1/6$, $p \ge 1$ and $T > 0$, there is a constant $K >0$ depending on $T$ such that
\beqnn
 \mbb{E}[|\nu_t^m(x) - \nu_r^m(x)|^{2p}]\le K(t - r)^{\delta + (p - 1)/2}, \quad 0 < r < t \le T.
\eeqnn
\end{lemma}
\proof
 By \eqref{num} we get
\beqnn
\ar\ar|\nu_t^m(x) - \nu^m_r(x)|\cr
\ar\ar\ =\Bigg| \int_0^r\int_{\mbb{R}}\sum_{k = 1}^m 1_{[t_{k - 1}, t_k)}(s) \sqrt{\gamma_m(\mu_{t_{k - 1}}^m, z)} G_m(\mu_s^m(z))\cr
\ar\ar\qquad\qquad\qquad\qquad\qquad\qquad \cdot[p_{r - s}(x-z) - p_{t- s}(x-z)] W(d s, d z)\cr
\ar\ar\ \quad - \int_r^t\int_{\mbb{R}}\sum_{k = 1}^m 1_{[t_{k - 1}, t_k)}(s) \sqrt{\gamma_m(\mu_{t_{k - 1}}^m, z)} G_m(\mu_s^m(z))p_{t\!  -\!  s}(x\!  -\!  z) W(d s, d z)\Bigg|.
\eeqnn
By the above, \eqref{gamma0722}, \eqref{G_m} and Burkholder-Davis-Gundy's inequality one can check that, for $0 < r < t \le T,$
\beqlb\label{0425b}
\ar\ar\mbb{E}\left[|\nu_t^m(x) - \nu^m_r(x)|^{2p}\right]\cr
\ar\ar\quad\le  K \mbb{E}\left[\left|\int_0^r\int_{\mbb{R}}\left[p_{r - s}(x - z) - p_{t - s}(x - z)\right]^2G_m(\mu_s^m(z))^2 d s d z\right|^p\right]\cr
\ar\ar \qquad+ K\mbb{E}\left[\Big|\int_r^t\int_{\mbb{R}}G_m(\mu_s^m(z))^2 p_{t - s}(x - z)^2d s d z\Big|^p\right]\cr
\ar\ar\quad\le  K \mbb{E}\left[\left|\int_0^r\int_{\mbb{R}}\left[p_{r - s}(x - z) - p_{t - s}(x - z)\right]^2[\mu_s^m(z) + 1] d s d z\right|^p\right]\cr
\ar\ar \qquad+ K\mbb{E}\left[\Big|\int_r^t\int_{\mbb{R}}[\mu_s^m(z) + 1] p_{t - s}(x - z)^2d s d z\Big|^p\right]\cr
\ar\ar\quad=: I_1^m(t,r) + I_2^m(t,r).
\eeqlb
Further, by H\"{o}lder's inequality,  Lemma \ref{l2.4}, \eqref{p_t} and \eqref{J},  we have
\beqlb\label{addd1}
\ar\ar\mbb{E}\left[\int_{\mbb{R}}p_{r - s}(x - z)^{2 - \delta} |\mu_s^m(z)|^{p} d z\right]\cr
\ar\ar\qquad \le K\mbb{E}\left[\left(\int_{\mbb{R}}|\mu_s^m(z)|^{2p}J(z) d z\right)^{1/2}\right] \left(\int_{\mbb{R}}p_{r - s}(x - z)^{4 - 2\delta}e^{|z|}d z\right)^{1/2}\cr
\ar\ar\qquad \le K (r - s)^{(2\delta - 3)/4}.
\eeqlb
Replacing $r$ with $t$ the above inequality still holds.
By H\"{o}lder's inequality,  \eqref{p_{r - s}} and \eqref{addd1}, for $\delta \in (0, 1/6)$ and $0 < r < t \le T,$ we have
\beqlb\label{0414}
\ar\ar\mbb{E}\left[\int_0^rd s\int_{\mbb{R}}[p_{r - s}(x - z) - p_{t - s}(x - z)]^2|\mu_s^m(z)|^{p}d z\right] \cr
\ar\ar\quad\le K(t - r)^{\delta}\int_0^r(r - s)^{-3\delta/2}d s\cr
\ar\ar\qquad\quad \cdot\mbb{E}\left[\int_{\mbb{R}}[p_{r - s}(x - z)^{2 - \delta} + p_{t - s}(x - z)^{2 - \delta}]|\mu_s^m(z)|^{p} d z\right]\cr
\ar\ar\quad\le K(t - r)^{\delta}\left[\int_0^r(r - s)^{-(3+4\delta)/4}d s + \int_0^r(r - s)^{-3\delta/2}(t - s)^{(2\delta - 3)/4}d s\right]\cr
\ar\ar\quad\le K(t - r)^{\delta}\left[T^{(1-4\delta)/4} \!+\! \left(\!\int_0^r(r - s)^{-6\delta}d s\!\right)^{1/4}\!\left(\!\int_0^r(t - s)^{(2\delta - 3)/3}d s\!\right)^{3/4}\!\right]\cr
\ar\ar\quad\le  K(t - r)^{\delta},
\eeqlb
where the above $K$ in the last inequality depends on $T.$ By \cite[Lemma~1.4.4]{X13a}, there exists a constant $K$ independent of $r, t, T$ such that
\beqlb\label{aaa}
\int_0^r\int_{\mbb{R}}[p_{r - s}(x - z) - p_{t - s}(x - z)]^2 d s d z \le K|t - r|^{1/2}
\eeqlb
and
\beqlb\label{a}
\int_r^t\int_{\mbb{R}} p_{t - s}(x - z)^2d s d z \le K|t - r|^{1/2}.
\eeqlb
 By H\"{o}lder's inequality, \eqref{0414} and \eqref{aaa}, it implies that
\beqlb\label{0425}
I_1^m(t,r) \ar\le\ar K\mbb{E}\left[\int_0^rd s\int_{\mbb{R}}[p_{r - s}(x - z) - p_{t - s}(x - z)]^2|\mu_s^m(z)|^{p}d z\right]\cr
\ar\ar\quad \cdot \left[\int_0^r\int_{\mbb{R}}[p_{r - s}(x - z) - p_{t - s}(x - z)]^2 d s d z\right]^{p - 1} + K(t - r)^{p/2}\ccr
\ar\le\ar K (t - r)^{\delta + (p-1)/2}  + K(t - r)^{p/2},
\eeqlb
where the above $K$ depends on $T.$ Similarly, by H\"{o}lder's inequality and Lemma~\ref{l2.4}, one can see that
\beqnn
\ar\ar\mbb{E}\left[\int_r^t \int_{\mbb{R}}|\mu_s^m(z)|^{p} p_{t - s}(x - z)^2 d s d z\right]\cr
\ar\ar\qquad \le K\mbb{E}\left[\int_r^t d s\left[\int_{\mbb{R}}|\mu_s^m(z)|^{2p} J(z)d z\right]^{1/2}\left[\int_{\mbb{R}}p_{t - s}(x - z)^4 e^{|z|} d z\right]^{1/2}\right]\cr
\ar\ar\qquad \le K\int_r^t (t - s)^{-3/4}d s = K(t - r)^{1/4},
\eeqnn
where the above $K$ depends on $T.$ By the above and \eqref{a} we obtain
\beqlb\label{0425a}
I_2^m(t,r) \ar\le\ar K\mbb{E}\left[\Big|\int_r^t\int_{\mbb{R}}\mu_s^m(z) p_{t - s}(x - z)^2d s d z\Big|^p\right] \cr
\ar\ar + K \Big|\int_r^t\int_{\mbb{R}} p_{t - s}(x - z)^2d s d z\Big|^p \cr
\ar\le\ar K \mbb{E}\left[\int_r^t\! \int_{\mbb{R}}|\mu_s^m(z)|^{p} p_{t - s}(x \!-\! z)^2 d s d z\right]\!\left[\int_r^t\! \int_{\mbb{R}} p_{t - s}(x \!-\! z)^2 d s d z\right]^{p - 1}\ccr
\ar\ar  + K(t - r)^{p/2}\ccr
\ar\le\ar K(t - r)^{1/4 + (p - 1)/2}  + K(t - r)^{p/2}.
\eeqlb
The result follows from \eqref{0425b}, \eqref{0425} and \eqref{0425a}.
\qed

Recall that $C([0, T],  C(\mbb{R})^+)$ is the space of the continuous maps from $[0, T]$ to $C(\mbb{R})^+$ with topology induced by \eqref{d1fg} by changing $\|f(t) - g(t)\|_{\mathcal{Y}}$ to $\|f(t) - g(t)\|_{C(\mbb{R})^+} = \int_0^\infty e^{-L}(\|f(t, \cdot) - g(t, \cdot)\|_{[-L, L]} \wedge 1)\d L,$ where $\|\cdot\|_{[-L, L]}$ is the supremum norm on $[-L, L]$. One can check that $C([0, T] \times \mbb{R}, \mbb{R}_+)$ is a subspace of $C([0, T],  C(\mbb{R})^+)$. Similar to above, we have the following result.

\begin{lemma}\label{l22.7}
For fixed $0 < \beta < 1$ and $p \ge 1$, there is a constant $K$ such that
\beqnn
\mbb{E}\left[|\nu_t^m(x) - \nu_t^m(y)|^{2p}\right]\le K|x - y|^{\beta p}, \qquad  x, y \in \mbb{R}.
\eeqnn
\end{lemma}

{\begin{lemma}\label{l0114}
Suppose that $(\mu_t^m)_{t \in [0, T]}$ is a solution to \eqref{mun1}. Then for every $m \ge 1$, we have $\mbb{E}[\<\mu_t^m, 1\>] = \<\mu_0, 1\>.$
\end{lemma}

\proof
Notice that $\gamma_m, G_m$ are bounded for every $m \ge 1.$ 
Let $t \in [0, T]$ be fixed. Recall that $p_{t-s}(x - z) = 0$ for all $s > t.$ Then for any $u \ge 0,$ we have 
\beqnn
\ar\ar\mbb{E}\left[\int_0^u\int_{\mbb{R}} \sum_{k = 1}^m1_{[t_{k-1}, t_k)}(s)\gamma_m(\mu_{t_{k - 1}}^m, z) G_m(\mu_s^m(z))^2p_{t - s}(x - z)^2 d s d z\right]\cr
\ar\ar\qquad\le K \int_0^u\int_{\mbb{R}} p_{t - s}(x - z)^2 d s d z \le K\sqrt{T}< \infty,
\eeqnn
which implies that 
\beqnn
u \mapsto \int_0^u\int_{\mbb{R}} \sum_{k = 1}^m1_{[t_{k-1}, t_k)}(s)\sqrt{\gamma_m(\mu_{t_{k - 1}}^m, z)} G_m(\mu_s^m(z))p_{t - s}(x - z) W(d s, d z) 
\eeqnn
is a martingale for any fixed $t \in [0, T]$. By \eqref{mun1}, one sees that $\mbb{E}[\mu_t^m(x)] = \<\mu_0, p_t(x - \cdot)\>$. Then $\mbb{E}[\<\mu_t^m, 1\>] = \<\mu_0, 1\>.$ The result follows.
\qed
}

\noindent{\it Proof of Theorem \ref{t1}.} By \eqref{mun1} and \eqref{num} one can check that
\beqnn
\mu_t^m(x) \ar=\ar \<\mu_0, p_t(x - \cdot)\> + \nu_t^m(x).
\eeqnn
Notice that $p$ can be taken greater than $\max\{2(1-\delta)+1, 1/\beta\}$ in Lemmas \ref{l22.6} and \ref{l22.7}. Then the Kolmogorov-Chentsov criterion (see \cite[Theorem 3.23]{K02}) holds. By \cite[Corollary 16.9]{K02}, one can check that, for each $T, L > 0,$ the sequence of $\{\mu_t^m(x): t \in [0, T], x\in [-L, L]\}$ in $C([0, T] \times [-L, L], \mbb{R}_+)$ is tight, which implies $\{\mu_t^m(x): t \in [0, T], x \in \mbb{R}\}$ is tight in $C([0, T] \times \mbb{R}, \mbb{R}_+)$ under the topology induced by \eqref{dfg}, hence, has a weakly convergent subsequence $\{\mu_t^{m_k}(x): t \in [0,T], x \in \mbb{R}\}$ with limit $\{\mu_t(x): t \in [0,T], x \in \mbb{R}\}$. Let
\beqnn
B_t(\phi) = \int_0^t \int_{\mbb{R}} \phi(x) W(d s, d x).
\eeqnn
Let $\mathcal{X}$ be a suitable Banach space containing $L^2(\mbb{R})$ such that $(\iota, L^2(\mbb{R}), \mathcal{X})$ being an abstract Wiener space, where $\iota$ denotes the inclusion map of $L^2(\mbb{R})$ into $\mathcal{X}$, see \cite[\S I.4]{K75} for the definition of the abstract Wiener space.  Similar to \cite[Theorem~3.2.4]{KX95}, one can check that $(B_t)_{t \ge 0}$ is an $L^2(\mbb{R})$-cylindrical Brownian motion taking values in $\mathcal{X}$. Thus, $(\mu^{m_k}, B) \rightarrow (\mu, B)$ in law on $C([0, T] \times \mbb{R}, \mbb{R}_+)\times C([0, T],  \mathcal{X})$ as $k \rightarrow \infty.$ Then, $(\mu_t^{m_k}, B_t)_{t \in [0,T]}\to (\mu_t, B_t)_{t \in [0,T]}$ in law on $C([0,T],C(\mbb{R})^+)\times C([0, T],  \mathcal{X})$.

 Applying Skorokhod's representation theorem (e.g. \cite[p.102, Theorem~1.8]{EK86}), on another probability space, there are continuous processes $(\hat{\mu}_t^{m_k}, \hat{B}_t^{m_k})_{t \in [0,T]}$ and $(\hat{\mu}_t, \hat{B}_t)_{t \in [0,T]}$
with the same distribution as  $(\mu_t^{m_k}, B_t)_{t \in [0,T]}$ and  $(\mu_t, B_t)_{t \in [0,T]},$ respectively. Moreover, $(\hat{\mu}_t^{m_k}, \hat{B}_t^{m_k})_{t \in [0,T]} \rightarrow (\hat{\mu}_t, \hat{B}_t)_{t \in [0,T]}$
 almost surely as $k \rightarrow \infty$.  Recall that there is a pathwise unique continuous $C_{\rm tem}^+(\mbb{R})$-valued solution to \eqref{mun}. For any $f \in C_c^\infty(\mbb{R})$, almost surely the following holds for each $t \in [0, T]$:
 \beqnn
 \<\hat{\mu}_t^{m_k}, f\> \ar=\ar \<\mu_0, f\> + \frac12\int_0^t \<\hat{\mu}_s^{m_k}, f''\>d s \\
 \ar\ar  + \sum_{j = 1}^{m_k}\int_{t_{j - 1}\wedge t}^{t_j\wedge t}\int_{\mbb{R}} \sqrt{\gamma_{m_k}(\hat{\mu}_{t_{j - 1}}^{m_k}, z)} G_{m_k}(\hat{\mu}_s^{m_k}(z))f(z) \hat{W}^{m_k}(d s, d z),
 \eeqnn
where $\hat{W}^{m_k}(d s, d z)$ is the corresponding time-space white noise of $(\hat{B}_t^{m_k})_{t \in [0,T]}$. It follows from Lemma~\ref{l2.4} and Fatou's lemma that, for every $p \ge 1,$ we have
\beqnn
 \sup_{0 \le t \le T}\mbb{E}\left[\int_{\mbb{R}}\hat{\mu}_t(x)^{2p}J(x)d x\right] \ar=\ar  \sup_{0 \le t \le T}\mbb{E}\left[\int_{\mbb{R}}\lim_{k \rightarrow \infty}\hat{\mu}_t^{m_k}(x)^{2p}J(x)d x\right]\cr
 \ar\le\ar \sup_{0 \le t \le T}\varliminf_{k \rightarrow \infty}\mbb{E}\left[\int_{\mbb{R}}\hat{\mu}_t^{m_k}(x)^{2p}J(x)d x\right]\cr
 \ar\le\ar \sup_{0 \le t \le T, k \ge 1}\mbb{E}\left[\int_{\mbb{R}}\hat{\mu}_t^{m_k}(x)^{2p}J(x)d x\right]\cr
 \ar\le\ar Ke^{KT} < \infty.
\eeqnn
Then \eqref{moment} holds for $(\hat{\mu}_t)_{t \ge 0}.$ Further, for every $t \in [0, T]$, by the above and Lemma \ref{l2.4} we have
\beqnn
\lim_{k \rightarrow \infty}\mbb{E}\left[\int_{\mbb{R}}|\hat{\mu}_t^{m_k}(x) - \hat{\mu}_t(x)|^{2}J(x)d x\right] = 0
\eeqnn
and
\beqnn
\sup_{0 \le t \le T, k \ge 1}\mbb{E}\left[\int_{\mbb{R}}[\hat{\mu}_t(x) + \hat{\mu}_t^{m_k}(x)]^{2}J(x)d x\right] < \infty.
\eeqnn
Then by the above we have
\beqnn
\ar\ar\mbb{E}\left[\int_0^t\langle|\hat{\mu}_s^{m_k} - \hat{\mu}_s|, f''\rangle d s\right]\cr
\ar \ar\quad \le { K}\mbb{E}\left[\int_0^td s\int_{\mbb{R}}|\hat{\mu}_s^{m_k}(x) - \hat{\mu}_s(x)|^2J(x)d x\right]^{1/2} \rightarrow 0
\eeqnn
as $k \rightarrow \infty.$ Recall that $\gamma$ satisfies Condition \ref{0.4b}. Thus,
\beqnn
\ar\ar \mbb{E}\left[\left|\int_0^t\int_{\mbb{R}} \sqrt{\gamma(\hat{\mu}_s, z)\hat{\mu}_s(z)}f(z) \hat{W}^{m_k}(d s, d z)\right|^2\right]\\
\ar\ar \qquad \le K\mbb{E}\left[\int_0^t\int_{\mbb{R}} \hat{\mu}_s(z)J(z) d s d z\right] \cr
\ar\ar \qquad \le K\left\{\mbb{E}\left[\int_0^t\int_{\mbb{R}} \hat{\mu}_s(z)^2J(z) d s d z\right]\right\}^{1/2} < \infty.
\eeqnn
By \cite[Lemma~2.4]{XY19}, for any $f \in C_c^\infty(\mbb{R})$, we have
\beqlb\label{0113}
 \<\hat{\mu}_t, f\> \ar=\ar \<\mu_0, f\> + \frac12\int_0^t \<\hat{\mu}_s, f''\>d s\cr
 \ar\ar + \int_0^t \int_{\mbb{R}} \sqrt{\gamma(\hat{\mu}_s, z)\hat{\mu}_s(z)}f(z) \hat{W}(d s, d z), \quad t \in [0, T]
\eeqlb
almost surely, where $\hat{W}(d s, d z)$ is the corresponding time-space white noise of $(\hat{B}_t)_{t \in [0, T]}$. { 
 By Fatou's lemma and Lemma \ref{l0114}, we have 
\beqnn
\mbb{E}[\<\hat{\mu}_t, 1\>] \ar=\ar \mbb{E}\left[\<\lim_{k \rightarrow \infty}\hat{\mu}_t^{m_k}, 1\>\right]\le \varliminf_{k \rightarrow \infty} \mbb{E}[\<\hat{\mu}_t^{m_k}, 1\>]\cr
 \ar=\ar \varliminf_{k \rightarrow \infty} \mbb{E}[\<\mu_t^{m_k}, 1\>] = \<\mu_0, 1\>.
\eeqnn
For any $n \ge m,$ let $1_{[m,n]}(x)$ be the characteristic function of $[m,n].$ By the above and dominated convergence theorem, we have 
\beqlb\label{0126a}
\lim_{m \rightarrow \infty}\lim_{n \rightarrow \infty}\mbb{E}[\<\hat{\mu}_t, 1_{[m,n]}\>] = 0.
\eeqlb 

For any fixed nonnegative function $f \in C_b^2(\mbb{R}),$ we will prove below that \eqref{0113} holds for any $t \in [0, T]$ almost surely. In fact, there exists a nonnegative function $f_{m,n} \in C_c^\infty(\mbb{R})$ such that 
\beqlb\label{0210}
\left\{
\begin{aligned}
&	f_{m,n}(x) =f(x), \qquad x \in [m,n];\\
&	0 \le f_{m,n}(x) \le f(x)\cdot 1_{[m-1, n+1]}(x);\\
&	|f_{m,n}''(x)| \le K1_{[m-1, n+1]}(x),
\end{aligned}\right.
\eeqlb
where $K > 0$ is a constant independent of $m, n.$ There is an example in $C_c^\infty(\mbb{R})$ satisfying \eqref{0210} below:
\beqnn
f_{m,n}(x) = \left\{
\begin{aligned}
	&0, \qquad\qquad\qquad\qquad\qquad\qquad\qquad\qquad\quad  x < m-1,\\
	&2f(x)\cdot\int_{m-1}^x \rho(2(y-m)+1) d y, \qquad\qquad  x \in [m-1, m),\\
	&f(x),\qquad\qquad\qquad\qquad\qquad\qquad\qquad\qquad  x \in [m,n),\\
	&f(x)\cdot\left[1-2\int_n^x \rho(2(y-n)-1)d y\right], \qquad  x \in [n, n+1),\\
	&0, \qquad\qquad\qquad\qquad\qquad\qquad\qquad\qquad\quad  x \ge n+1,
\end{aligned}\right.
\eeqnn
where $\rho$ is the mollifier defined by \eqref{rho}. It is easy to see that $0 \le f_{m,n}(x) \le f(x)\cdot 1_{[m-1, n+1]}(x).$
Moreover, $\sup_{x \in (-1, 1)}\rho(x) \le C$ and $\sup_{x \in (-1, 1)} |\rho'(x)|\le 8C,$ where $C$ is the constant in \eqref{rho}. For any $x \in [m-1,m),$ 
\beqnn
|f_{m,n}''(x)| \ar=\ar \Big|2f''(x)\cdot\int_{m-1}^x \rho(2(y-m)+1) d y + 4f'(x)\cdot\rho(2(x-m)+1) \cr
\ar\ar + 4f(x)\cdot\rho'(2(x-m)+1)\Big|\cr
\ar\le\ar 2\|f''\| + 4C\|f'\| + 32C\|f\|=:K,
\eeqnn
where $\|\cdot\|$ is the supremum norm.  Similarly, the above also holds for any $x \in [n,n+1)$. This means $ |f_{m,n}''(x)| \le K1_{[m-1, n+1]}(x).$ One sees that the above $f_{m,n}$ satisfies \eqref{0210}. By \eqref{0126a}, \eqref{0210} and the fact of $f(x) \le \|f\|$ for all $x \in \mbb{R}$, we have
\beqlb\label{0207}
\lim_{m \rightarrow \infty}\lim_{n \rightarrow \infty}\mbb{E}\left[\<\hat{\mu}_t, f_{m,n}\>+\<\hat{\mu}_t, |f_{m,n}''|\>\right] = 0.
\eeqlb
Moreover, it follows from Doob's inequality that
\beqlb\label{eq0206a}
\ar\ar\mbb{E}\left[\sup_{0 \le t \le T}\left|\int_0^t \int_{\mbb{R}} \sqrt{\hat{\mu}_s(z)\gamma(\hat{\mu}_s, z)}f_{m,n}(z)\hat{W}(d s, d z)\right|^2\right]\cr
\ar\ar\qquad\le \mbb{E}\left[\int_0^T \int_{\mbb{R}} \hat{\mu}_s(z)\gamma(\hat{\mu}_s, z)f_{m,n}(z)^2d s d z\right]\cr
\ar\ar\qquad\le K\int_0^T\mbb{E}[ \<\hat{\mu}_s, f_{m,n}\>]d s .
\eeqlb
By \eqref{0113}, it holds that
\beqnn
\mbb{E}\left[\sup_{0\le t\le T}\<\hat{\mu}_t, f_{m,n}\>\right] \ar\le\ar \<\mu_0, f_{m,n}\> +\frac{1}{2} \mbb{E}\left[\int_0^T\<\hat{\mu}_s, |f_{m,n}''|\>d s\right]\cr
\ar\ar + \mbb{E}\left[\sup_{0 \le t \le T}\left|\int_0^t \int_{\mbb{R}} \sqrt{\hat{\mu}_s(z)\gamma(\hat{\mu}_s, z)}f_{m,n}(z)\hat{W}(d s, d z)\right|\right].
\eeqnn
Recall that $\mu_0 \in C_c(\mbb{R})^+.$  Taking $n \rightarrow \infty$ and then $m \rightarrow \infty,$ by \eqref{0207}, \eqref{eq0206a} and  the above we have
\beqlb\label{eq0206b}
\lim_{m\rightarrow\infty}\lim_{n\rightarrow\infty}\mbb{E}\left[\sup_{0\le t\le T}\<\hat{\mu}_t, f_{m,n}\>\right]=0.
\eeqlb
On the other hand, by Fatou's lemma and \eqref{0210} we have
\beqnn
\<\hat{\mu}_t, f1_{[m,\infty)}\> \le \varliminf_{n\rightarrow\infty}\<\hat{\mu}_t, f1_{[m,n]}\> \le \varliminf_{n\rightarrow\infty}\<\hat{\mu}_t, f_{m,n}\>,
\eeqnn
which implies that
\beqnn
\mbb{E}\left[\sup_{0\le t\le T}\<\hat{\mu}_t, f1_{[m,\infty)}\>\right] \ar\le\ar \mbb{E}\left[\sup_{0\le t\le T}\varliminf_{n\rightarrow\infty}\<\hat{\mu}_t, f_{m,n}\>\right]\cr
\ar\le\ar \varliminf_{n\rightarrow\infty}\mbb{E}\left[\sup_{0\le t\le T}\<\hat{\mu}_t, f_{m,n}\>\right].
\eeqnn
By the above and \eqref{eq0206b}, one can see that
\beqlb\label{0206}
\lim_{m\rightarrow \infty} \mbb{E}\left[\sup_{0\le t\le T}\<\hat{\mu}_t, f\cdot 1_{[m,\infty)}\>\right] = 0.
\eeqlb
Similarly, 
\beqlb\label{0206a}
\lim_{m\rightarrow \infty}\mbb{E}\left[\sup_{0\le t\le T}\<\hat{\mu}_t, f\cdot 1_{(-\infty, -m]}\>\right] = 0.
\eeqlb
Now we take the sequence of functions $\{f_{-m,m}\}_{m=1}^\infty$ satisfying \eqref{0210}. It follows from \eqref{0210}, \eqref{0206} and \eqref{0206a} that
\beqnn
\mbb{E}\left[\sup_{0 \le t \le T}|\<\hat{\mu}_t, f_{-m,m}-f\>|\right] \!\!\ar\le\ar\!\! \mbb{E}\left[\sup_{0 \le t \le T}\int_{\mbb{R}}\hat{\mu}_t(x)|f_{-m,m}(x) - f(x)|d x\right]\cr
	\ar\le\ar \!\!\mbb{E}\left[\sup_{0 \le t \le T}\<\hat{\mu}_t, f\cdot1_{[m,\infty)}+ f\cdot1_{(-\infty, -m]}\>\right]
	\eeqnn
	goes to $0$ as $m\rightarrow\infty$. Moreover, by \eqref{0210} and dominated convergence theorem, we have 
	\beqnn
	\lim_{m\rightarrow \infty}\!\!\mbb{E}\left[\left|\int_0^T\!\! \<\hat{\mu}_s, f_{-m,m}'' - f''\>d s\right|\right]\!\! \ar\le\ar \!\!	\lim_{m\rightarrow \infty}\mbb{E}\!\!\left[\int_0^T \!\!\<\hat{\mu}_s, |f_{-m,m}'' - f''|\>d s\right]\cr
	\ar=\ar\!\! \mbb{E}\left[\int_0^T \!\!\!\<\hat{\mu}_s, \lim_{m\rightarrow \infty}\!\!|f_{-m,m}'' - f''|\>d s\right]\! = \!0.
	\eeqnn
	Similarly, one can check that $\lim_{m\rightarrow \infty} [|\<\mu_0, f_{-m,m} - f\>|] =0$ and
\beqnn
\ar\ar\mbb{E}\left[\sup_{0 \le t \le T}\left(\int_0^t \int_{\mbb{R}} \sqrt{\gamma(\hat{\mu}_s, z)\hat{\mu}_s(z)}(f_{-m,m}(z) -f(z)) \hat{W}(d s, d z)\right)^2\right]\cr
\ar\ar\qquad \le K \mbb{E}\left[\int_0^T \int_{\mbb{R}}\hat{\mu}_s(z)(f_{-m,m}(z) -f(z))^2d s d z \right]\rightarrow 0
\eeqnn
as $m \rightarrow \infty.$ Then for any nonnegative function $f \in C_b^2(\mbb{R}),$ \eqref{0113} holds for any $t \in [0, T]$ almost surely. 

For any $f \in C_b^2(\mbb{R}),$ there exist nonegative functions $f_1, f_2 \in C_b^2(\mbb{R})$ such that $f =f_1-f_2.$
 Then for any $f \in C_b^2(\mbb{R}),$ \eqref{0113} holds for any $t \in [0, T]$ almost surely. Letting $T \rightarrow \infty$ in \eqref{0113},} 
it completes the proof of the existence of the solution to \eqref{1.1}. By Proposition~\ref{t2.4} one can see that $X_t(d x) = \mu_t(x)d x$ satisfies the MP (\ref{0.1}, \ref{0.2}), which implies the conclusion.
\qed

\noindent{\it Proof of Theorem~\ref{t2}.}
Assume that $(\mu_t)_{t \ge 0}$ is a solution to SPDE \eqref{1.1} satisfying \eqref{moment}. Then it also satisfies \eqref{mu} by \cite[Theorem~2.1]{S94}. Let $r_0 > 0$ be fixed. For any $r_0 \le r \le t \le T$ and $p \ge 1$, we have
\beqlb\label{I123}
\ar\ar\mbb{E}\big[|\mu_t(x) - \mu_r(x)|^{2p}\big]\cr
\ar \ar\quad\le K\Big|\int_\mbb{R} [p_t(x - z)- p_r(x - z)]\mu_0(z) d z\Big|^{2p}\cr
\ar\ar\qquad + K\mbb{E}\left[\Big|\int_r^t\int_{\mbb{R}}p_{t - s}(x - z)\sqrt{\mu_s(z)}W(d s, d z)\Big|^{2p}\right]\cr
\ar\ar\qquad+ K\mbb{E}\left[\Big|\int_0^r \int_{\mbb{R}}[p_{r - s}(x - z) - p_{t - s}(x - z)]\sqrt{\mu_s(z)}W(d s, d z)\Big|^{2p}\right]\cr
\ar\ar\quad=: K(I_1 + I_2 + I_3).
\eeqlb
 Recall that $\mu_0 \in C_c(\mbb{R})^+$ and then $\int_{\mbb{R}}\mu_0(z)^2d z < \infty.$ Note that \eqref{p_{r - s}} holds for $\delta = 1$ by \cite[Lemma~III 4.5]{P02}.
	By \eqref{p_t}, \eqref{p_{r - s}} and
H{\"o}lder's inequality, we have
\beqlb\label{I1a}
I_1 \ar=\ar \Big|\int_\mbb{R} [p_t(x - z)- p_r(x - z)]\mu_0(z)d z\Big|^{2p}\cr
\ar\le \ar \Big|\int_\mbb{R} [p_t(x - z)- p_r(x - z)]^2 d z\int_{\mbb{R}}\mu_0(z)^2d z\Big|^{p}\cr
\ar\le\ar K \Big|\int_\mbb{R} [p_t(x - z)- p_r(x - z)]^2 d z\Big|^{p}\cr
\ar\le\ar K(t - r)^{ p}r^{-3 p/2}\Big|\int_\mbb{R} \left[p_t(x - z) + p_r(x - z)\right] d z\Big|^p\cr
\ar\le \ar K(t - r)^{p}r_0^{-3 p/2}, \qquad r_0 \le r \le t \le T.
\eeqlb
The following estimations of $I_2$ and $I_3$ are similar with that of $I_1^m(t, r)$ and $I_2^m(t, r)$ in Lemma \ref{l22.6}, respectively. By \eqref{moment}, \eqref{hhhh} and H{\"o}lder's inequality, for $\delta \in (0, 1/6)$ one can see that
\beqlb\label{kk}
\ar\ar\mbb{E}\left[\int_{\mbb{R}}p_{t - s}(x - z)^{2-\delta} \mu_s(z)^p d z \right]\cr
\ar\ar\qquad \le \mbb{E}\left[\left(\int_{\mbb{R}}|\mu_s(z)|^{2p}e^{-|z|}d z\right)^{1/2}\right]\left[\int_{\mbb{R}}p_{t - s}(x - z)^{4-2\delta}e^{|z|}d z\right]^{1/2}\cr
\ar\ar\qquad \le K (t - s)^{(2\delta - 3)/4}
\eeqlb
and
\beqnn
\ar\ar\mbb{E}\left[\int_r^t\int_{\mbb{R}}p_{t - s}(x - z)^2 \mu_s(z)^pd s d z \right]\cr
\ar\ar\quad \le \mbb{E}\left[\int_r^t d s \left[\int_{\mbb{R}}|\mu_s(z)|^{2p}e^{-|z|}d z\right]^{1/2}\left[\int_{\mbb{R}}p_{t - s}(x - z)^{4}e^{|z|}d z\right]^{1/2}\right]\cr
\ar\ar\quad \le K\int_r^t (t - s)^{-3/4}d s \le K(t - r)^{1/4},
\eeqnn
which implies that
\beqlb\label{I2a}
I_2 \ar\le\ar K \mbb{E}\left[\Big|\int_r^t\int_{\mbb{R}}p_{t - s}(x - z)^2 \mu_s(z)d s d z\Big|^p\right]\cr
\ar\le\ar K \mbb{E}\left[\int_r^t\int_{\mbb{R}}p_{t - s}(x - z)^2 \mu_s(z)^pd s d z \right]\left[\int_r^t\int_{\mbb{R}}p_{t - s}(x - z)^2d s d z\right]^{p - 1}\cr
\ar\le\ar K(t - r)^{1/4 + (p - 1)/2},
\eeqlb
where the above constant $K$ only depends on $T.$  By \eqref{p_{r - s}}, \eqref{kk} and H\"{o}lder's inequality, for $\delta \in (0, 1/6)$ we have
\beqnn
\ar\ar\mbb{E}\left[\int_0^rd s\int_{\mbb{R}}[p_{r - s}(x - z) - p_{t - s}(x - z)]^2\mu_s(z)^{p}d z\right] \cr
\ar\ar\quad\le K(t - r)^{\delta}\int_0^r(r - s)^{-3\delta/2}d s\cr
\ar\ar\qquad\quad \cdot\mbb{E}\left[\int_{\mbb{R}}[p_{r - s}(x - z)^{2 - \delta} + p_{t - s}(x - z)^{2 - \delta}]\mu_s(z)^{p} d z\right]\cr
\ar\ar\quad\le K(t - r)^{\delta}\left[\int_0^r(r - s)^{-(3+4\delta)/4}d s + \int_0^r(r - s)^{-3\delta/2}(t - s)^{(2\delta - 3)/4}d s\right]\cr
\ar\ar\quad\le K(t - r)^{\delta}\left[T^{(1-4\delta)/4} + \left(\int_0^r(r - s)^{-6\delta}d s\right)^{1/4}\left(\int_0^r(t - s)^{(2\delta - 3)/3}d s\right)^{3/4}\right]\cr
\ar\ar\quad\le  K(t - r)^{\delta}
\eeqnn
with $K$ only depending on $T.$ Combining the above with \eqref{aaa}, it implies that
\beqlb\label{I3a}
I_3 \ar\le\ar K \mbb{E}\left[\Big|\int_0^r\int_{\mbb{R}}[p_{r - s}(x - z) - p_{t - s}(x - z)]^2\mu_s(z) d s d z\Big|^p\right]\cr
\ar\le\ar K \mbb{E}\left[\int_0^r\int_{\mbb{R}}[p_{r - s}(x - z) - p_{t - s}(x - z)]^2\mu_s(z)^p d s d z\right]\cr
\ar\ar\qquad\cdot \left[\int_0^r\int_{\mbb{R}}[p_{r - s}(x - z) - p_{t - s}(x - z)]^2 d s d z\right]^{p - 1}\cr
\ar\le\ar K (t - r)^{\delta + (p -1)/2}
\eeqlb
with $\delta \in (0, 1/6)$ and the constant $K$ only depending on $T.$ By \eqref{I123}, \eqref{I1a}, \eqref{I2a} and \eqref{I3a}, we have
\beqlb\label{tr}
\mbb{E}\big[|\mu_t(x) - \mu_r(x)|^{2p}\big]
\le K(t - r)^{\delta + (p - 1)/2},
\eeqlb
where $\delta \in (0, 1/6)$ and the above constant $K$ only depends on $r_0, T.$ Similarly,
\beqlb\label{xy}
\mbb{E}\left[|\mu_t(x) - \mu_t(y)|^{2p}\right] \le K|x - y|^{\beta p}, \quad  \forall\ x, y \in \mbb{R}
\eeqlb
with $\beta \in (0, 1),$ where the above constant $K$ depends on $r_0.$ By \eqref{tr} and \eqref{xy}, for $t, r \in [r_0, T]$ and $x, y \in \mbb{R},$ it is easy to see that
\beqnn
\mbb{E}\left[|\mu_t(x) - \mu_r(y)|^{2p}\right] \ar\le\ar K\mbb{E}\left[|\mu_t(x) - \mu_t(x)|^{2p}\right] + K\mbb{E}\left[|\mu_t(x) - \mu_t(y)|^{2p}\right]\cr
\ar\le\ar K(t - r)^{\delta + (p - 1)/2} + K|x - y|^{\beta p}
\eeqnn
Taking $p$ large enough, the result follows from the Kolmogorov's continuity criteria (see, e.g., \cite[Theorem~3.23]{K02}).
\qed

\section{Proof of Theorem~\ref{main}}

In this section we show the weak uniqueness of the solution to the MP (\ref{0.1}, \ref{0.2}) under Condition \ref{0.4b}. The main idea is to relate the MP~(\ref{0.1}, \ref{0.2}) with a system of SPDEs, which is satisfied by a sequence of corresponding distribution-function-valued processes. The weak uniqueness of the solution to the MP (\ref{0.1}, \ref{0.2}) follows from the pathwise uniqueness of the solution to the SPDEs.

 For $(X_t)_{t \ge 0}$ satisfying the MP (\ref{0.1}, \ref{0.2}), we define the $[a_i, a_{i + 1})$-distribution-function-valued process
\beqlb\label{u}
 u_t^i(x) = X_t((a_i, x]), \qquad  a_i \le x < a_{i + 1}, \ i = 0, \cdots, n.
\eeqlb
In Proposition~\ref{t1.2}, for $i = 0, \cdots, n$, we show the processes $(u_t^i)_{t \ge 0}$ defined by \eqref{u} solves the following system of SPDEs:
\begin{numcases}{}
\begin{split}\label{2.1}
u_t^i(x) &= u_0^i(x) + \int_0^t \frac{1}{2}\Delta u_s^i(x) d s + \int_0^t\int_0^{u_s^i(x)}g_i(\nabla u^i_s(a_{i + 1}))W_i(d s, d z),\\
&\qquad\qquad \qquad\qquad \qquad \qquad\quad x \in [a_i, a_{i + 1}),\ i = 0, \cdots, n - 1;
\end{split}
\\
\begin{split}\label{2.1a}
u_t^n(x) &= u_0^n(x) + \int_0^t \frac{1}{2}\Delta u_s^n(x) d s + \int_0^t\int_0^{u_s^n(x)}g_n(u_s^n(\infty))W_n(d s, d z),\\
& \qquad\qquad \qquad\qquad \qquad \qquad\qquad \qquad\qquad\qquad\quad x \in [a_n, \infty).
\end{split}
\\
\begin{split}\label{0.7}
&u^i(a_i)=0,\;\; \;i=0,1,\cdots, n,\\
&\nabla u^{i-1}_t(a_i)=\nabla u^i_t(a_i),\;\;\;i=1,\cdots, n.
\end{split}
\end{numcases}
 where $u_s^n(\infty):=\lim_{x \rightarrow \infty}u_s^n(x)$ and $W_i(d s, d z),\ i = 0, \cdots, n$ are independent time-space Gaussian white noises on $\mbb{R}_+ \times \mbb{R}_+$ with intensity $d s d z$. The pathwise uniqueness of the solution to the SPDEs is obtained in Proposition \ref{t1.3}.

The system of SPDEs \eqref{2.1}-\eqref{0.7} will be understood in the following form: for any $\phi_i \in C_b^2[a_i, a_{i + 1}]$ with $\phi_i(a_i) = \phi'_i(a_{i + 1}) = 0,\ i = 0, 1, \cdots, n - 1,$
and $\phi_n \in C_b^2[a_n, \infty)$ with $\phi_n(a_n) = \phi_n(\infty) = 0$ (given $\phi_n(\infty) := \lim_{x \rightarrow \infty}\phi_n(x)$),  almost surely for each $t \ge 0,$ we have
{ \begin{numcases}{}
\begin{split}\label{u^i}
\langle u_t^i, \phi_i\rangle &= \langle u_0^i, \phi_i\rangle + \frac{1}{2}\int_0^t \big[\langle u_s^i, \phi_i''\rangle + \phi_i(a_{i + 1})\nabla u^i_s(a_{i + 1})\big]d s\\
 & + \int_0^t\int_0^\infty\left[\int_{a_i}^{a_{i + 1}} 1_{\{z \le u^i_s(x)\}}\phi_i(x)d x\right]g_i(\nabla u^i_s(a_{i + 1}))W_i(d s, d z);
\end{split}
\\
\begin{split}\label{un}
\langle u_t^n, \phi_n\rangle &= \langle u_0^n, \phi_n\rangle + \frac{1}{2}\int_0^t \langle u_s^n, \phi_n''\rangle d s\\
 &  + \int_0^t\int_0^\infty\left[\int_{a_n}^\infty 1_{ \{z \le u^n_s(x)\}}\phi_n(x) d x\right]g_n\big(u_s^n(\infty)\big)W_n(d s, d z),
\end{split}
\end{numcases}}
where $\langle u_0^i, \phi_i\rangle = \int_{\mbb{R}}u_0^i(x)\phi_i(x)d x$ and $\langle u_0^n, \phi_n\rangle = \int_{\mbb{R}}u_0^n(x)\phi_n(x)d x.$

\begin{prop}\label{t1.2}
Suppose that $(X_t)_{t \ge 0}$ is a solution to the MP (\ref{0.1}, \ref{0.2}). Then $\{u_t^i(x): t \ge 0, x \in [a_i, a_{i + 1})\}, i = 0, 1, \cdots, n$ defined as \eqref{u}, solves the group of SPDEs (\ref{2.1},\ref{2.1a}) with boundary condition \eqref{0.7}.
\end{prop}
\proof
 For any $\phi_i\in C^3_c(a_i,a_{i+1})$ with $i = 0, 1, \cdots, n-1$, by integration by parts, almost surely for each $t \ge 0$ we have
 \beqlb\label{0330}
\<u^i_t,\phi_i'\> \ar=\ar - \<X_t,\phi_i\>=-M_t(\phi_i) - \<X_0,\phi_i\> - \frac12\int^t_0\<X_s,\phi_i''\>ds\cr
\ar=\ar - M_t(\phi_i) + \<u^i_0,\phi_i'\> + \frac12\int^t_0\<u^i_s, \phi_i'''\> ds.
\eeqlb
Thus
\beqlb\label{eq0325b}
-M_t(\phi_i)=\<u^i_t,\phi_i'\> -\<u^i_0,\phi_i'\> -\frac12\int^t_0\<u^i_s, (\phi_i')''\> ds
\eeqlb
is a continuous martingale. By Lemma~\ref{t2.3} we have
\beqnn
- M_t(\phi_i) = \int_0^t\int_{\mbb{R}}\phi_i(x)M(d s, d x)
\eeqnn
 with quadratic variation process $(\<-M(\phi_i)\>_t)_{t \ge 0}$ defined by
\beqnn
\<-M(\phi_i)\>_t\ar=\ar \int^t_0g_i(\nabla u^i_s(a_{i + 1}))^2ds\int_{\mbb{R}} \phi_i(x)^2 X_s(d x)\cr
\ar=\ar \int^t_0g_i(\nabla u^i_s(a_{i + 1}))^2d s \int_{0}^{u_s^i(a_{i + 1})}\phi_i(u_s^i(y)^{-1})^2 d y\cr
\ar=\ar \int^t_0\int^\infty_0\(\int_{a_i}^{a_{i + 1}} 1_{\{y\le u^i_s(x)\}}\phi_i'(x)dx\)^2g_i(\nabla u^i_s(a_{i + 1}))^2ds dy,
\eeqnn
where $u_s^i(y)^{-1}$ denotes the generalized inverse of the nondecreasing function $u_s^i,$ that is, $u_s^i(y)^{-1} = \inf\{x \in [a_i, a_{i + 1}): u_s^i(x) \ge y\}.$
Moreover, for $\phi_n \in C_c^3(a_n, \infty)$, one can see that
\begin{equation}\label{eq0325n}
-M_t(\phi_n) = \<u^n_t,\phi_n'\> -\<u^n_0,\phi_n'\> -\frac12\int^t_0\<u^n_s, \phi_n'''\> ds
\end{equation}
is a continuous martingale with quadratic variation process $(\<-M(\phi_n)\>_t)_{t \ge 0}$ defined by
\beqnn
\<-M(\phi_n)\>_t\ar=\ar \int^t_0g_n(u_s^n(\infty))^2ds\int_{\mbb{R}} \phi_n(x)^2 X_s(d x)\\
\ar=\ar \int^t_0\int^{\infty}_0\(\int_{a_n}^{\infty} 1_{\{y \le u^n_s(x)\}}\phi_n'(x)dx\)^2g_n(u_s^n(\infty))^2ds dy.
\eeqnn
As in Lemma~\ref{t2.3}, the family $\{-M_t(\phi_i): t \ge 0, \phi_i \in C_c^3(a_i, a_{i + 1})\}$ determines a martingale measure $\{M_t(B): t \ge 0, B \in \mcr{B}(a_i, a_{i + 1}), i = 0, \cdots, n-1\}.$ Moreover, for $\phi_i \in C_c^3(a_i, a_{i + 1})$ and $\phi_j \in C_c^3(a_j, a_{j + 1})$ with $i \neq j,$ we have $\phi_i(x)\phi_j(x) = 0$ for any $i \neq j,$ and then by Lemma \ref{t2.3},
\beqnn
\ar\ar\<-M(\phi_i), -M(\phi_j)\>_t\cr
\ar \ar\quad=\int_0^t d s \int_{\mbb{R}} \gamma(\mu_s, z)X_s(d z)\int_{\mbb{R}}\phi_i(x)\delta_z(d x)  \int_{\mbb{R}}\phi_j(y)\delta_z(d y)\\
\ar\ar\quad= \int_0^t d s \int_{\mbb{R}} \gamma(\mu_s, z)\phi_i(z)\phi_j(z)X_s(d z) = 0.
\eeqnn
By El Karoui and M\'{e}l\'{e}ard \cite[Theorem~III-7, Corollary~III-8]{EM90}, on some extension of the probability space, one can define a sequence independent Gaussian white noises $W_i(d s, d z), i = 0, \cdots, n$ on $(0, \infty)^2$ based on $d s d z$ such that, for any $\phi_i \in C_c^3(a_i, a_{i + 1}), i = 0, \cdots, n - 1$,
\beqlb\label{03301}
-M_t(\phi_i)\ar=\ar \int^t_0\int^\infty_0\left[\int_{a_i}^{a_{i + 1}} 1_{\{z \le u^i_s(x)\}}\phi_i'(x) \sqrt{\gamma(\mu_s, x)}dx\right]W_i(d s, d z)\cr
\ar=\ar\int_{a_i}^{a_{i + 1}}\phi_i'(x) dx \int^t_0\int^{u^i_s(x)}_0g_i(\nabla u^i_s(a_{i + 1}))W_i(d s, d z)
\eeqlb
for each $t \ge 0$ almost surely, and for any  $\phi_n \in C_c^3(a_n, \infty)$, we have
\beqlb\label{03302}
-M_t(\phi_n) = \int_{a_n}^{\infty} \phi_n'(x) dx\int^t_0\int^{u^n_s(x)}_0g_n( u^n_s(\infty))W_n(d s, d z)
\eeqlb
for each $t \ge 0$ almost surely. By \eqref{0330}, \eqref{03301} and \eqref{03302}, the equations \eqref{u^i}, \eqref{un} hold for any $\phi_i \in C_c^2(a_i, a_{i + 1}),\ i = 0, \cdots, n.$ Moreover, for each $T > 0,$
\beqnn
\sup_{0 \le t\le T}\sup_{x \in [a_i, a_{i +1}]}|u_t^i(x)| = \sup_{0 \le t \le T} X_t([a_i, a_{i+1}]) < \infty
\eeqnn
almost surely for $i = 0, 1, \cdots, n-1.$ Therefore, for any $\phi_i \in C_b^2[a_i, a_{i + 1}]$ with $\phi_i(a_i) = \phi_i'(a_{i+1}) = 0,\ i = 0, \cdots, n - 1,$ \eqref{u^i} follows from Lemma~\ref{l3.6} in the Appendix. Moreover, for any $\phi_n \in C_b^2[a_n, \infty)$ with $\phi_n(a_n) = \phi_n(\infty) = 0,$ \eqref{un} follow from Lemma~\ref{l0330} in the Appendix.
\qed

\begin{lemma}\label{l3.2}
 Suppose that $(u_t^n(x))_{t \ge 0, x \ge a_n}$ satisfies \eqref{2.1a}  with $u_t^n(a_n) = 0$. Let $u_t^n(\infty) = \lim_{x \rightarrow \infty}u_t^n(x).$ Then $(u_t^n(\infty))_{t \ge 0}$ satisfies
 \beqlb\label{2.10}
u_t^n(\infty)=u_0^n(\infty)+\int_0^t
\int_0^{u_s^n(\infty)}g_n(u_s^n(\infty))W_n(d s, d z).
 \eeqlb
\end{lemma}

\proof
Recall that $p_t(x)=\frac{1}{\sqrt{2\pi}t}\e^{-x^2/(2t)}$ and let
\beqnn
q_t^x(y):=p_t(x+a_n-y)-p_t(x-a_n+y)
 \eeqnn
for $t>0$ and $x,y\ge a_n$. One can check that $q_t^x(a_n) = q_t^x(\infty) = 0$ for any $t > 0$ and $x \ge a_n.$
Then \eqref{2.1a} can be written in the following mild form, whose proof is similar to \cite[Lemma~5.1]{XY20}, or \cite[Theorem~7.26]{L11} taking into consideration of the boundary condition:
\beqlb\label{2.8}
u_t^n(x)\ar=\ar \int_0^t\int_0^\infty
\Big[\int_{a_n}^\infty1_{\{z \le u_s^n(y)\}}q_{t-s}^x(y) d y\Big]g_n(u_s^n(\infty)) W_n(d s, d z)\cr
\ar\ar + \<u_0^n,q_t^x\>
 \eeqlb
for any $x \ge a_n,$ where $\<u_0^n,q_t^x\> = \int_{\mbb{R}}u_0^n(y)q_t^x(y)d y.$ By a change of variable and dominated convergence theorem, we have
 \beqnn
\<u_0^n,q_t^x\>
 \ar=\ar
\int_{a_n}^\infty u_0^n(y)[p_t(x+a_n-y)-p_t(x-a_n+y)]d y \\
 \ar=\ar
\int_{-\infty}^x u_0^n(x+a_n-z)p_t(z)d z
-\int_x^\infty u_0^n(z-x+a_n)p_t(z)d z \\
 \ar\to\ar
u_0^n(\infty)\int_{-\infty}^\infty p_t(z)d z=u_0^n(\infty)
 \eeqnn
and
\beqnn
\ar\ar\int_{a_n}^\infty1_{\{z \le u_s^n(y)\}}q_{t-s}^x(y) d y \cr
\ar\ar\qquad= \int_{a_n}^\infty 1_{\{z \le u_s^n(y)\}}[p_t(x+a_n-y)-p_t(x-a_n+y)]d y \cr
\ar\ar\qquad= \int_{-\infty}^x 1_{\{z \le u_s^n(x + a_n - y)\}}p_t(y)d y - \int_x^\infty 1_{\{z \le u_s^n(y - x + a_n)\}}p_t(y) d y\cr
\ar\ar\qquad \to 1_{\{z \le u_s^n(\infty)\}}\int_{-\infty}^\infty p_t(y)d y = 1_{\{z \le u_s^n(\infty)\}}
 \eeqnn
as $x\to\infty,$ which ends the proof.
\qed

\begin{prop}\label{t1.3}(Pathwise uniqueness)
 Suppose that \eqref{beta} holds, and $(u_t)_{t \ge 0}$, $(\tilde{u}_t)_{t \ge 0}$ are two solutions with the same initial value to (\ref{2.1}, \ref{2.1a}) satisfying the boundary conditions~\eqref{0.7}. Then $\mbb{P}\{u_t(x) = \tilde{u}_t(x)\ \text{for all}\ t \ge 0,\ x \in \mbb{R}\} = 1.$
\end{prop}

\proof
By Lemma~\ref{l3.2}, $(u_t^n(\infty))_{t \ge 0}$ and $(\tilde{u}_t^n(\infty))_{t \ge 0}$ are two solutions to \eqref{2.10} with the same initial value. Then by \eqref{beta} and \cite[Theorem 2.1]{DL12}, the pathwise uniqueness of the solution to \eqref{2.10} holds, i.e., $\mbb{P}\{u_t^n(\infty) = \tilde{u}_t^n(\infty) \ \text{for all}\ t \ge 0\} = 1.$
Further, the pathwise uniqueness of the solution holds for \eqref{2.1a} by Lemma~\ref{lemma4.3}, i.e., $\mathbb{P}\{u_t^n(x) = \tilde{u}_t^n(x)\  \text{for\ all}\ t \ge 0,\ x \ge a_n\} = 1.$ That implies the strong uniqueness of $(\nabla u_t^n(a_n))_{t \ge 0}.$ By Lemma~\ref{lemma4.4} and using induction, one can get the pathwise uniqueness of the solution to \eqref{2.1} for $i = 0, 1, \cdots, n-1.$ The proof ends here.
 \qed

\noindent{\it Proof of Theorem~\ref{main}.}
	Suppose that $(X_t)_{t \ge 0}$ and $(\tilde{X}_t)_{t \ge 0}$ are two solutions to the MP (\ref{0.1}, \ref{0.2}). In fact, the existence of the solution to the MP follows from Theorem~\ref{t1}. As \eqref{u}, one can define $u_t^i(x)$ and $\tilde{u}_t^i(x)$ corresponding to $X_t$ and $\tilde{X}_t,$ respectively. In fact, by Proposition \ref{t1.2} one can see that $(u_t^i(x))_{t \ge 0, x \in [a_i, a_{i +1})}$ and $(\tilde{u}_t^i(x))_{t \ge 0, x \in [a_i, a_{i +1})},\ i = 0, 1, \cdots, n$ are two solutions to the group of SPDEs (\ref{2.1}, \ref{2.1a}) with boundary condition \eqref{0.7}. Further, the strong uniqueness of the solution to (\ref{2.1}, \ref{2.1a}) follows from Proposition~\ref{t1.3}, which implies the weak uniqueness of the solution to the MP (\ref{0.1}, \ref{0.2}).
\qed

\section{Appendix}
In this section we give some useful lemmas which are used in the proofs of Propositions \ref{t1.2} and \ref{t1.3} . 
Let $\Phi\in C_c^2(0,1)$ satisfy $0\le \Phi\le 2$ and $\int_{0}^1\Phi(x) d x=1$. 
For $k\ge1$ and $x\in[0,1]$ let
 \beqlb\label{hn}
h_k(x):= \int_0^{kx}\Phi(z) d z\cdot \int_{x^k}^1\Phi(z) d z.
 \eeqlb
Then $h_k\in C_c^2(0,1)$ for all $k\ge1$.

\begin{lemma}\label{l3.5}
Suppose that $f\in C[0,1]$ with $f'(1)$ and $f'(0)$ existing.
Then
 \beqnn
\lim_{k\to\infty}\<f,h_k'\> =  f(0)-f(1),\quad
\lim_{k\to\infty}\<f,h_k''\> = f'(1)-f'(0),
 \eeqnn
where $\<f, g\> = \int_0^1 f(x)g(x)d x$ with $f, g\in C[0, 1].$
\end{lemma}
\proof
Observe that for each $n\ge1$,
 \beqnn
h_k'(x)=k\Phi(kx)
\int_{x^k}^1\Phi(z) d z
-kx^{k-1}\Phi(x^k)\int_0^{kx}\Phi(z)d z,\qquad x\in[0,1].
 \eeqnn
Then by change of variables and dominated convergence,
as $k\to\infty$ we have
 {\small\beqnn
\ar\ar\<f,h_k'\>\cr
 \ar\ar\quad=\int_0^1f(x)k\Phi(kx)\Big[\int_{x^k}^1\Phi(z) d z\Big] d x
-\int_0^1f(x)kx^{k-1}\Phi(x^k)\Big[\int_0^{kx}\Phi(z) d z\Big] d x \cr
 \ar\ar\quad=
\int_0^1f(y/k)\Phi(y)\Big[\int_{(y/k)^k}^1\Phi(z) d z\Big] d y
-\int_0^1f(y^{1/k})\Phi(y)\Big[\int_0^{ky^{1/k}}\Phi(z) d z\Big] d y
 \eeqnn}
converges to $f(0)-f(1)$ as $k \rightarrow \infty$,
which gives the first assertion.

In the following we prove the second assertion.
Observe that
 \beqlb\label{5.1}
h_k''(x)
 \ar=\ar
k^2\Phi'(kx)\int_{x^k}^1\Phi(z) d z-2k^2 x^{k-1}\Phi(kx)\Phi(x^k) \cr
\ar\ar
-\big[k(k-1)x^{k-2}\Phi(x^k)
+k^2x^{2k-2}\Phi'(x^k)\big]\int_0^{kx}\Phi(z) d z \cr
 \ar=:\ar
M_{1,k}(x)-2M_{2,k}(x) - M_{3,k}(x).
 \eeqlb
By change of variables and dominated convergence again,
as $k\to\infty$, we have
 \beqlb\label{5.2}
\ar\ar\int_0^1[f(x)-f(0)]M_{1,k}(x) d x\cr
 \ar\ar\quad=
\int_0^1k[f(y/k)-f(0)]\Phi'(y)\Big[\int_{y^kk^{-k}}^1\Phi(z) d z\Big] d y  \cr
 \ar\ar\quad\to
f'(0)\int_0^1y \Phi'(y) d y
=-f'(0)
 \eeqlb
and
 \beqlb\label{5.3}
\ar\ar\int_0^1[f(x)-f(1)]M_{2,k}(x) d x\cr
 \ar\ar\quad  =\int_0^1\frac{f(y^{1/k})-f(1)}{y^{1/k}-1}k(y^{1/k}-1)\Phi(ky^{1/k})\Phi(y) d y\cr
\ar\ar\quad\to 0.
 \eeqlb
Similarly, as $k\to\infty$,
 \beqlb\label{5.4}
 &\ &
\int_0^1[f(x)-f(1)]M_{3,k}(x) d x \cr
 &\ &\quad=
\int_0^1\frac{f(y^{1/k})-f(1)}{y^{1/k}-1}k(y^{1/k}-1)
\big[k^{-1}(k-1)y^{-1/k}\Phi(y) \cr
 &\ &\quad\qquad
+y^{(k-1)/k}\Phi'(y)\big]\cdot\Big[\int_0^{ky^{1/k}}\Phi(z) d z\Big] d y\cr
 &\ &\quad\to
f'(1)\int_0^1 \ln y[\Phi(y)+y\Phi'(y)] d y=-f'(1).
 \eeqlb
Applying integration by parts and the fact $0\le\Phi\le 2$ with $\mbox{supp}(\Phi)\subset(0,1)$, one can check that
 \beqnn
\int_0^1M_{1,k}(x) d x
 \ar=\ar
\int_0^1\Big(\int_0^{kx}\Phi(z) d z\Big)''\cdot \Big(\int_{x^k}^1\Phi(z) d z\Big) d x \\
 \ar=\ar \int_0^1M_{2,k}(x) d x =
k\int_0^1 \Phi(ky^{1/k})\Phi(y) d y\cr
\ar=\ar
k\int_0^{k^{-k}} \Phi(ky^{1/k})\Phi(y) d y
\le 4 k^{1-k}
 \eeqnn
and
 \beqnn
\int_0^1M_{3,k}(x) d x \ar=\ar
-\int_0^1\Big(\int_{x^k}^1\Phi(z) d z\Big)''\cdot\Big(\int_0^{kx}\Phi(z) d z\Big) d x\\
\ar=\ar \int_0^1M_{2,k}(x) d x
\le 4 k^{1-k}.
 \eeqnn
Then combining \eqref{5.1} with \eqref{5.2}-\eqref{5.4}
one completes the proof.
\qed

\begin{lemma}\label{l3.6}
Suppose that for each $\phi \in C_c^2(0, 1)$,  $(u_t)_{t\ge0}$ satisfies
 \beqlb\label{uiphi}
\<u_t,\phi\>
\ar=\ar
\<u_0,\phi\>  + \frac12\int_0^t\<u_s,\phi''\> d s\nonumber\\
&\ & +\int_0^t\int_0^\infty g (\nabla u_s(1)) \Big[\int_0^1 1_{\{z \le u_s (x)\}}\phi(x)d x\Big]W (d s, d z),
 \eeqlb
and for each $T > 0,$
\beqlb\label{u0130}
\sup_{0 \le t \le T}\sup_{x \in [0, 1]}|u_t(x)| < \infty \qquad a.s.
\eeqlb
Then for each $\phi\in C_b^2[0, 1]$, almost surely for $t \ge 0$ we have
 \beqnn
\<u_t,\phi\>
\ar=\ar
\<u_0,\phi\>  + \frac12\int_0^t\big[\<u_s,\phi''\>
+F_s(\phi)] d s\nonumber\\
&\ & + \int_0^t\int_0^\infty g(\nabla u_s(1))\Big[\int_0^11_{\{z \le u_s(x)\}}\phi(x) d x\Big]W(d s, d z),
 \eeqnn
 where
 \beqnn
F_s(\phi):=[\phi(1)\nabla u_s(1)-\phi(0)\nabla u_s(0)]-[u_s(1)\phi'(1)- u_s(0)\phi'(0)].
 \eeqnn
\end{lemma}

\proof
Recall $h_k$ in \eqref{hn}.
For $m\ge1$ define stopping time $\tau_m$ by
 \beqnn
\tau_m:=\inf\Big\{t\ge0:\sup_{x\in[0,1]}|u_t(x)|\ge m\Big\}
 \eeqnn
with the convention $\inf\emptyset=\infty$.
Then $\lim_{m \to\infty}\tau_m=\infty$
 almost surely by \eqref{u0130}.
It follows from \eqref{uiphi} that,  for any $\phi \in C_b^2[0, 1]$ and $m, k > 1,$
\beqlb\label{3.3}
\<u_{t\wedge\tau_m},\phi h_{k}\>
 \ar=\ar \int_0^{t\wedge\tau_m}\int_0^\infty g(\nabla u_s(1))\Big[\int_0^11_{\{z \le u_s(x)\}}\phi(x)h_{k}(x)d x\Big]W(d s, d z)\nonumber\\
 &\ &
+\<u_0,\phi h_{k}\> +\frac12\int_0^{t\wedge\tau_m}\<u_s,(\phi h_{k})''\> d s .
 \eeqlb
Notice that
 \beqnn
\<u_s,(\phi h_{k})''\>
=
\<u_s,\phi''h_{k}\>
+2\<u_s,\phi'h'_{k}\>
+\<u_s,\phi h''_{k}\>.
 \eeqnn
It follows from Lemma \ref{l3.5} that
 \beqnn
\ar\ar\lim_{k\to\infty}\<u_s,(\phi h_{k})''\>\cr
 \ar\ar\quad=
\<u_s,\phi''\>
+[u_s(0)\phi'(0)-u_s(1)\phi'(1)]
-[\phi(0)\nabla u_s(0)-\phi(1)\nabla u_s(1)] \\
 \ar\ar\quad=
\<u_s,\phi''\> + F_s(\phi).
 \eeqnn
Thus letting $k\to\infty$ in \eqref{3.3} we obtain
 \beqnn
\<u_{t\wedge\tau_m},\phi\>
\ar=\ar
\<u_0,\phi\>
+\frac12\int_0^{t\wedge\tau_m}\big[\<u_s,\phi''\>
+ F_s(\phi)\big]d s\\
&\ & +\int_0^{t\wedge\tau_m}\int_0^\infty g(\nabla u_s(1))\Big[\int_0^1 1_{\{z \le u_s(x)\}}\phi(x)d x\Big]W(d s, d z).
 \eeqnn
Letting $m\to\infty$ the result holds.
\qed

\begin{lemma}\label{l0330}
Suppose that for each $\phi \in C_c^2(0, \infty)$,  $(u_t)_{t\ge0}$ satisfies
 \beqnn
\<u_t,\phi\>
\ar=\ar
\<u_0,\phi\>  + \frac12\int_0^t\<u_s,\phi''\> d s\cr
&\ & +\int_0^t\int_0^\infty g (u_s(\infty)) \Big[\int_0^\infty 1_{\{z \le u_s (x)\}}\phi(x)d x\Big]W (d s, d z),
 \eeqnn
 and for each $T > 0,$
\beqlb\label{u0infty30}
\sup_{0 \le t \le T}\sup_{x \in [0, \infty)}|u_t(x)| < \infty \qquad a.s.
\eeqlb
Then for each $\phi\in C_b^2[0, \infty)$ with $\phi(\infty) = 0$, almost surely for $t \ge 0$ we have
 \beqnn
\<u_t,\phi\>
\ar=\ar
\<u_0,\phi\>  + \frac12\int_0^t\big[\<u_s,\phi''\>
+u_s(0)\phi'(0) - \phi(0)\nabla u_s(0)] d s\cr
&\ & + \int_0^t\int_0^\infty g(\nabla u_s(\infty))\Big[\int_0^\infty 1_{\{z \le u_s(x)\}}\phi(x) d x\Big]W(d s, d z).
 \eeqnn
\end{lemma}

\proof
We take $h_k(x) := \int_0^{kx}\Phi(z) d z,$ where $\Phi \in C_c^2(0, \infty)$, $\int_0^\infty \Phi(x)d x = 1$ and $0 \le \Phi \le 2.$ Simliar with the proof of Lemma~\ref{l3.5}, it is easy to check that
\beqnn
\<f, h_k'\> \rightarrow f(0) \qquad\qquad \<f, h_k''\> \rightarrow -f'(0)
\eeqnn
as $k \rightarrow \infty,$ where $f \in C_b[0, \infty)$ with $f'(0)$ existing and $f(\infty) = \lim_{x \rightarrow \infty}f(x) = 0.$ The following proof is similar with that of Lemma~\ref{l3.6}, we omit it here.
\qed

Next we consider the following equation:
\begin{numcases}{}
\begin{split}\label{a2}
u_t(x) &= u_0(x) + \int_0^t \frac{1}{2}\Delta u_s(x) d s + \int_0^t\int_0^{u_s(x)}G_sW(d s, d z), \quad x \ge 0,\\
u_t(0)&=0,
\end{split}
\end{numcases}
where $(G_t)_{t \ge 0}$ is a bounded continuous process, and $W(d s, d z)$ is a time-space Gaussian white noise on $\mbb{R}_+\times\mbb{R}_+$ with intensity $d s d z$. Further, we assume \eqref{u0infty30} holds for the above equation. The $[0, \infty)$-distribution-function-valued process $(u_t)_{t \ge 0}$ solving \eqref{a2} will be understood in the following form: for each $\phi \in C_b^2[0, \infty)$ with $\phi(0) = \phi(\infty) = 0$, we have
\beqnn
\<u_t, \phi\> \ar=\ar \<u_0, \phi\> + \frac{1}{2}\int_0^t\<u_s, \phi''\>d s\\
\ar\ar + \int_0^t\int_0^\infty\left[\int_0^\infty 1_{\{z \le u_s(x)\}}\phi(x) d x\right] G_s W(d s, d z),\quad t \ge 0
\eeqnn
almost surely. It is easy to see that $t \goto \<u_t, \phi\>$ is continuous almost surely for each $\phi \in C_b^2[0, \infty).$ The pathwise uniqueness of the solution to the general version of \eqref{a2} is considered in Xiong and Yang~\cite[Theorem 1.9]{XY20}. We give a brief proof below, since that paper has not been formally published yet.

\begin{lemma}\label{lemma4.3}
Let $(u_t^i)_{t \ge 0}, i = 1, 2$ be two $[0, \infty)$-distribution-function-valued solutions to \eqref{a2} with the same initial value satisfying \eqref{u0infty30}. Then $\mathbb{P}\{u_t^1(x) = u_t^2(x)\  \text{for\ all}\ t, x \ge 0\} = 1.$
\end{lemma}

\proof
For each $k \ge 1$ we define $a_k = \exp\{-k(k + 1)/2\}.$ Then $a_k \rightarrow \infty$ and $\int_{a_k}^{a_{k - 1}}z^{-1} d z = k.$ Let $x \rightarrow \psi_k(x)$ be a positive continuous function supported by $(a_k, a_{k - 1})$ such that $\int_{a_k}^{a_{k - 1}}\psi_k(x) d x = 1$ and $\psi_k(x) \le (2kx)^{-1}$ for every $x > 0.$ For $k \ge 1$ let
\beqnn
\phi_k(z) = \int_0^{|z|} d y\int_0^y\psi_k(x) d x, \quad z \in \mbb{R}.
\eeqnn
It is easy to see that $|\phi_k'(z)| \le 1$ and $0 \le |z|\phi_k''(z) = |z|\psi_k(|z|) \le 2k^{-1}$ for $z \in \mbb{R}.$ Moreover, we have $\phi_k(z) \rightarrow |z|$ increasingly as $k \rightarrow \infty.$

Recall that $(u_t^1)_{t \ge 0}$ and $(u_t^2)_{t \ge 0}$ are two solutions to \eqref{a2} with the same initial value. Let $q_t^x(y) = p_t(x, y) - p_t(-x, y).$ It follows that for $i = 1, 2,$
\beqlb\label{u_t^i}
\<u_t^i, q_\delta^x\> \ar=\ar \<u_0, q_\delta^x\> + \frac{1}{2}\int_0^t \Delta_x(\<u_s^i, q_\delta^x\>) d s\ccr
\ar\ar + \int_0^t\int_0^\infty \left[\int_0^\infty 1_{\{z \le u_s^i(y)\}}q_\delta^x(y)d y\right]G_s W(d s, d z),
\eeqlb
where $\Delta_x$ is the second order spatial differential operator with respect to the variable $x$. For $n \ge 1,$ we define the stopping time
\beqnn
\tau_n := \inf\{t \ge 0: \sup_{x \ge 0}|u_t^1(x)| + \sup_{x \ge 0}|u_t^2(x)| \ge n\}.
\eeqnn
Then $\tau_n \rightarrow \infty$ almost surely as $n \rightarrow \infty$ by \eqref{u0infty30}. Let $v_t(x) = u_t^1(x) - u_t^2(x)$ and $v_t^\delta(x) = \<v_t, q_\delta^x\>.$ From \eqref{u_t^i} it follows that
\beqnn
v_t^\delta(x) = \frac{1}{2}\int_0^t \Delta_x (v_s^\delta(x)) d s + \int_0^t\int_0^\infty M_s^\delta(x, z) W(d s, d z),
\eeqnn
where
\beqnn
M_s^\delta(x, z) = G_s \int_0^\infty (1_{\{z \le u_s^1(y)\}} - 1_{\{z \le u_s^2(y)\}})q_\delta^x(y)d y.
\eeqnn
It then follows from It\^{o}'s formula that
\beqnn
\phi_k(v_{t\wedge \tau_n}^\delta(x)) \ar=\ar \frac{1}{2}\int_0^{t\wedge \tau_n} \phi_k'(v_s^\delta(x)) \Delta_x (v_s^\delta(x)) d s\\
\ar\ar + \int_0^{t\wedge \tau_n}\int_0^\infty \phi_k''(v_s^\delta(x)) |M_s^\delta(x, z)|^2 d s d z\\
\ar\ar + \int_0^{t\wedge \tau_n}\int_0^\infty \phi_k'(v_s^\delta(x)) M_s^\delta(x, z) W(d s, d z).
\eeqnn
Recall that $J(x) = \int_{\mbb{R}}e^{-|y|}\rho(x - y)d y$ satisfying \eqref{J}.  It leads to
\beqlb\label{delta}
\ar\ar\int_{\mbb{R}} J(x)\mathbb{E}[\phi_k(v_{t\wedge \tau_n}^\delta(x))] d x\cr
\ar\ar\quad = \frac{1}{2}\mbb{E}\left[\int_0^{t\wedge \tau_n} \int_{\mbb{R}} J(x)\phi_k'(v_s^\delta(x)) \Delta_x (v_s^\delta(x)) d xd s\right]\cr
\ar\ar\quad\quad  +\mbb{E}\left[ \int_0^{t\wedge \tau_n}\int_0^\infty\int_{\mbb{R}} J(x) \phi_k''(v_s^\delta(x)) |M_s^\delta(x, z)|^2 d s d zd x\right]\cr
\ar\ar\quad =: \mbb{E}\left[\int_0^{t\wedge \tau_n} \Big(\frac{1}{2}I_{1,k}^\delta(s) +  I_{2,k}^\delta(s) \Big)d s\right].
\eeqlb
By integration by parts and \eqref{J}, we have
\beqnn
I_{1,k}^\delta(s) \ar=\ar \int_{\mbb{R}} \Delta_x(\phi_k(v_s^\delta(x)))J(x)d x - \int_{\mbb{R}}  \phi_k''(v_s^\delta(x))|\nabla_xv_s^\delta(x)|^2J(x)d x\\
\ar\le\ar  \int_{\mbb{R}} \Delta_x(\phi_k(v_s^\delta(x))J(x)d x  =  \int_{\mbb{R}} \phi_k(v_s^\delta(x))J''(x)d x\\
\ar\le\ar K  \int_{\mbb{R}} \phi_k(v_s^\delta(x))J(x)d x.
\eeqnn
Recall that $(G_t)_{t \ge 0}$ is bounded. Then we have
\beqnn
I_{2,k}^\delta(s)
\le K \int_\mbb{R} J(x) \phi_k''(v_s^\delta(x))\<|v_s^\delta|, q_\delta^x\> d x.
\eeqnn
Letting $\delta \rightarrow 0$ in \eqref{delta} and using dominated convergence we have
\beqnn
\ar\ar\int_\mbb{R} J(x)\mathbb{E}[\phi_k(v_{t\wedge \tau_n}(x))] d x\\
 \ar\ar\quad\le K \int_0^td s \int_\mbb{R} \mbb{E} [\phi_k(v_{s\wedge \tau_n}(x))]J(x)d x\\
\ar\ar\quad\quad + K\int_0^td s \int_\mbb{R} \mbb{E} \big[\phi_k''(v_{s\wedge \tau_n}(x))|v_{s\wedge \tau_n}(x)|\big]J(x)d x.
\eeqnn
Recall that $0 \le |y|\phi_k''(y) \le 2/k$ for all $y \in \mbb{R},$ letting $k \rightarrow \infty$ in the above inequality we get
\beqnn
\int_\mbb{R} J(x)\mathbb{E}[|v_{t\wedge \tau_n}(x)|] d x \le K \int_0^t d s \int_\mbb{R} J(x)\mathbb{E}[|v_{s\wedge \tau_n}(x)|] d x,
\eeqnn
which implies that $\int_\mbb{R} J(x)\mathbb{E}[|v_{t\wedge \tau_n}(x)|] d x = 0$. Letting $n \rightarrow \infty,$ we obtain that $\int_\mbb{R} J(x)\mathbb{E}[|v_{t}(x)|] d x = 0$ by Fatou's lemma. Note that $t \rightarrow \int_\mbb{R} J(x)u_t^i(x) d x$ is continuous almost surely for $i = 1, 2$. Then $\mbb{P}\{\int_\mbb{R} J(x)|v_{t}(x)| d x = 0,\ \text{for\ all}\ t \ge 0\} = 1.$ The result follows since $(u_t^1)_{t \ge 0}$ and $(u_t^2)_{t \ge 0}$ are two $[0, \infty)$-dsitribution-function-valued processes.
\qed

Now we consider the following equation:
\begin{numcases}{}
\begin{split} \label{u01}
u_t(x) &= u_0(x) + \int_0^t \frac{1}{2}\Delta u_s(x) d s  +\int_0^t\int_0^{u_s(x)}\tilde{G}_s W(d s, d z), \quad x \in [0, 1],\\
u_t(0) &= 0, \quad \nabla u_t(1) = \mu_t,
\end{split}
\end{numcases}
where $(\mu_t)_{t \ge 0}$ and $(\tilde{G}_t)_{t \ge 0}$ are continuous processes, and $(\tilde{G}_t)_{t \ge 0}$ is bounded. Moreover, $W(d s, d z)$ is a time-space Gaussian white noise on $\mbb{R}_+\times\mbb{R}_+$ with intensity $d s d z$. The $[0, 1]$-distribution-function-valued process $(u_t(x))_{t \ge 0, x \in [0, 1]}$ solving \eqref{u01} will be understood in the following form:  for each $\phi \in C_b^2[0, 1]$ with $\phi(0) = \phi'(1) = 0,$ we have
\beqnn
\<u_t, \phi\> \ar=\ar \<u_0, \phi\> + \frac{1}{2}\int_0^t \<u_s, \phi''\> + \phi(1)\mu_s d s\\
\ar\ar + \int_0^t\int_0^\infty\left[\int_0^\infty 1_{\{z \le u_s(x)\}}\phi(x) d x\right] \tilde{G}_s W(d s, d z), \quad t \ge 0
\eeqnn
almost surely. The pathwise uniqueness of the solution to the general version of \eqref{u01} is considered in Xiong and Yang~\cite[Theorem 1.4]{XY20}.

\begin{lemma}\label{lemma4.4}
 Let $(u_t^i(x))_{t \ge 0, x \in [0,1]}, i = 1, 2$ be two $[0,1]$-distribution-function-valued solutions to \eqref{u01} with the same initial value satisfying \eqref{u0130}. Then
 \beqnn
 \mathbb{P}\left\{u_t^1(x) = u_t^2(x)\  \text{for  all} \ t\ge 0, x \in [0, 1]\right\} = 1.
 \eeqnn
\end{lemma}
\proof
The conclusion can be justified by using essentially the same argument as that in the proof of Lemma~\ref{lemma4.3} with $q_t^x(y)$ replaced by
\beqnn
q_t^x(y) \ar=\ar 2\sum_{k = -\infty}^\infty [p_t(4k + x - y)- p_t(4k - x - y)]\cr
\ar\ar - \sum_{k = - \infty}^{\infty} [p_t(2k + x - y) - p_t(2k - x - y)].
\eeqnn
We omit the details. \qed

\medskip

{\bf Acknowledgements} { We would like to express our sincere gratitude to the editors and an anonymous referee for their very helpful comments on the paper.} The research of L. Ji was supported in part by Guangdong Basic
and Applied Basic Research Foundation 2021A1515010031; the research of J. Xiong  was supported in part by NSFC grants 61873325, 11831010 and SUSTech fund Y01286120; the research of X. Yang was supported in part by NSFC grants 12061004, NSF of Ningxia grant 2021AAC02018 and the Major Research Project for North Minzu University No. ZDZX201902. 

\end{document}